\documentclass[a4paper,11pt]{article}

 \usepackage{graphicx}
\usepackage[width=17cm,height=25cm]{geometry}
 \usepackage{epsfig}
\usepackage{amssymb}
\usepackage{amsmath}

\newcommand{\be}{\begin{equation}}
\newcommand{\ee}{\end{equation}}
\newcommand{\bea}{\begin{eqnarray}}
\newcommand{\eea}{\end{eqnarray}}
\newcommand{\bbea}{\begin{eqnarray*}}
\newcommand{\eeea}{\end{eqnarray*}}

\newtheorem{theorem}{Theorem}[section]
\newtheorem{lemma}{Lemma}[section]

\newtheorem{example}{Example}[section]
\newtheorem{definition}{Definition}[section]
\newtheorem{proposition}{Proposition}[section]
\newtheorem{corollary}{Corollary}[section]

\newcommand{\gy}{\psi}

\begin{document}

\title{Spectral theory of a class of Block Jacobi matrices  and applications}

\author{Jaouad Sahbani}
\date{
Institut de Math\'ematiques de Jussieu-Paris Rive Gauche-UMR7586\\
 Universit{\'e} Paris Diderot (Paris 7) \\
 B\^atiment Sophie Germain--case 7012\\
5 rue Thomas Mann, 75205 Paris Cedex 13, FRANCE \\
{jaouad.sahbani@imj-prg.fr}}

\maketitle
\begin{abstract}
We develop a spectral analysis of a class of  block Jacobi operators based on the conjugate operator method of Mourre. 
We give several applications including scalar Jacobi operators
with periodic coefficients, a class of difference operators on cylindrical domains such as discrete wave propagators,  
and certain fourth-order difference operators.
\end{abstract}



%\tableofcontents

%%%%%%%%%%%%%%%%%%%%%%%%%%%%%%%%%%%%%%%%%%%%%%%%%%%%%%%%%%%%%%%%%%%%%%%%%%%%%%%%%%
%%%%%%%%%%%%%%%%%%%%%%%%%%%%%%%%%%%%%%%%%%%%%%%%%%%%%%%%%%%%%%%%%%%%%%%%%%%%%%%%%%
\section{Introduction}
%%%%%%%%%%%%%%%%%%%%%%%%%%%%%%%%%%%%%%%%%%%%%%%%%%%%%%%%%%%%%%%%%%%%%%%%%%%%%%%%%%
%%%%%%%%%%%%%%%%%%%%%%%%%%%%%%%%%%%%%%%%%%%%%%%%%%%%%%%%%%%%%%%%%%%%%%%%%%%%%%%%%%
%\subsection{Overview}
%%%%%%%%%%%%%%%%%%%%%%%%%%%%%%%%%%%%%%%%%%%%%%%%%%%%%
%%%%%%%%%%%%%%%%%%%%%%%%%%%%%%%%%%%%%%%%%%%%%%%%%%%%%
Let  ${\mathcal H}=l^2(\mathbb{Z}, \mathbb{C}^N)$, for some integer  $N\geq1$, be the Hilbert
space of square summable vector-valued sequences $ (\gy_n)_{n\in\mathbb{Z}} $ endowed
with the scalar product
$$
\langle\phi,\psi\rangle=\sum_{n\in\mathbb{Z}}\langle\phi_n,\psi_n\rangle_{\mathbb{C}^N}
 $$
 where $\langle\cdot,\cdot\rangle_{\mathbb{C}^N}$ is the usual scalar product of $\mathbb{C}^N$.
 Let $A_n$ and $B_n$ be two  sequences of $N$ by $N$ matrices such that $B_n=B_n^*$ for all $n\in\mathbb{Z}$. 
 Here we denote by $T^*$ the adjoint matrix of a given matrix $T$.
 Let us consider  the difference operator $H=H(\{A_n\},\{B_n\})$ acting in $\mathcal{H}$  by
\begin{equation}\label{Jacobi1}
(H\psi)_n=A_{n-1}^*\psi_{n-1}+B_n\psi_n+A_{n}\psi_{n+1}, ~~\mbox{
for all } n\in\mathbb{Z}.
\end{equation}
In the literature,  $H$ is usually called the block Jacobi operator realized by the sequences $A_n$ and $B_n$. 
Operators of such a form  are naturally related to matrix orthogonal polynomials theory, see for example \cite{B,C1,C,DPS,DeDe,K1,K2} 
and references therein.
Of course, the case where $N=1$ corresponds to the usual scalar Jacobi operators 
that are well studied by different approaches, see e.g. 
 \cite{B,BS1,BS2,CH,C,DKS,DK,D,DP,H,JN1,KvM,vM,P,S4,T} and their references. (This list is very restrictive 
 and contains only some papers related directly to the present one). Applications given in the second part of this paper represent
 additional motivations of our interest to block Jacobi operators.

Our main goal here is to study spectral properties of the
operator $H$ in terms of the asymptotic behavior of the
entries $A_n$ and $B_n$.  The key step of our analysis is the
 construction of a conjugate operator for $H$ in the sense of  the Mourre estimate, see \cite{M}.  
Recall that such an estimate has many important consequences on the spectral measure of $H$ as well as
 on the asymptotic behavior  of the resolvent of $(H-z)^{-1}$ when the complex parameter $z$ approaches the real axis, see 
 \cite{ABG,BG,BGS2,M,S1}.  

We need the following standard notations. 
For a self-adjoint operator $T$ we denote by $E_T(\cdot)$ its spectral measure, $\sigma(T)$ its spectrum,
$\sigma_{ess}(T)$ its essential spectrum,
 $\sigma_p(T)$ the set of its eigenvalues,  $\sigma_{sc}(T)$ its singular continuous spectrum, and
  $\sigma_{ac}(T)$ its absolutely continuous spectrum.

In this paper we will  focus on the so-called generalized Nevai class.
  More specifically,  we assume that
there exist two matrices $A$ and $B$ such that,
\begin{equation}\label{Jacobi2}
\lim_{|n|\rightarrow\infty}(\|A_n-A\|+\|B_n-B\|)=0.
\end{equation}
In particular, $H$ is a bounded self-adjoint operator in ${\mathcal H}$.  
Notice that in the literature, it is customary to assume that the matrices $A_n$ are positive definite or at least non singular.
Here we do not demand such an assumption, which turns out to be quite useful in many applications, see Section \ref{periodique}.

 It is more convenient to decompose $H$ into the sum $H=H_0+V$ where
 the operator $H_0=H_0(A,B)$ is defined   on $\mathcal{H}$ by
\begin{equation}\label{Jacobi0}
(H_0\psi)_n=A^*\psi_{n-1}+B\psi_n+A\psi_{n+1}, ~~\mbox{
for all } n\in\mathbb{Z}
\end{equation}
and  $V=H(\{A_n-A\},\{B_n-B\})$ is the Block Jacobi operator realized by the 
sequences $A_n-A,B_n-B$.
The operator $H_0$ is  clearly bounded  and self-adjoint  in $\mathcal{H}$ while,
  according  to (\ref{Jacobi2}),  the perturbation $V$ is  compact.
Hence, by Weyl theorem, $H$ and $H_0$ have the same essential spectra,
$$
\sigma_{ess}(H)=\sigma_{ess}(H_0).
$$
The point here is that $H_0$ is a block Toeplitz matrix whose  spectrum can be determined explicitly.  Indeed,  
by using Fourier transform we easily show  that $H_0$ is unitarily equivalent  to the direct integral 
\begin{equation}\label{direct}
\hat{H_0}=\int_{[-\pi,\pi]}^\oplus h(p)dp\quad\mbox{ acting in }\quad
\hat{\mathcal{H}}=\int_{[-\pi,\pi]}^\oplus\mathbb{C}^Ndp
\end{equation}
where the reduced operators $h(p)$ are given by,
\begin{equation}\label{reduced}
h(p)=e^{-ip}A+e^{ip}A^*+B.
\end{equation}
The  family  $h(p),p\in[-\pi,\pi]$, is  an  analytic family of self-adjoint matrices.
So, see \cite{Ka},  there exist two analytic
functions $\lambda_{j}$ and $W_{j}$ such that, for all $p\in[-\pi,\pi]$, 
$\{\lambda_{j}(p);\,1\leq j\leq N\}$ are the repeated eigenvalues  of
$h(p)$ and $\{W_{j}(p);\,1\leq j\leq N\}$ is a corresponding orthonormal basis
of eigenvectors. In particular, 
for any $j=1,\cdots,N$,  the critical set of $\lambda_j$  defined by 
$$
\kappa(\lambda_j)=\{\lambda_j(p)~|~ \lambda_j'(p)=0\}
$$
 is clearly  finite.
Therefore, the critical set of $H_0$ defined by
$$
\kappa(H_0)=\cup_{i=1}^{i=N}\kappa(\lambda_j)
$$
is finite too.  
We have
%
%
%x
\begin{proposition}\label{prop1}
The  spectrum of $H_0$ is purely absolutely continuous outside the critical set $\kappa(H_0)$. More precisely,
$
\sigma(H_0)=\cup_{i=1}^{i=N}\lambda_i([-\pi,\pi])~Ê,\quad\sigma_p(H_0)\subset\kappa(H_0)\quad\mbox{ and }~~\sigma_{sc}(H_0)=\emptyset.
$
\end{proposition}
%
%
%
%%%%%%%%%%%%%%%%%%%%%%%%%%%%%%%%%%%%%
%
Let $\alpha_j=\min_{p\in[-\pi,\pi]}\lambda_j(p)$ and 
$\beta_j=\max_{p\in[-\pi,\pi]}\lambda_j(p)$. Assume that  the $\lambda_j$'s are arranged 
so that $\alpha_1\leq\alpha_2\cdots\leq\alpha_N$.
So each $\lambda_j$ gives arise to a spectral  band
$\Sigma_j=\lambda_j([-\pi,\pi])=[\alpha_j,\beta_j]$,  $j=1,\cdots,N$. Two successive bands $\Sigma_j,\Sigma_{j+1}$ 
may overlap if $\alpha_{j+1}\leq\beta_j$; 
or may be disjoint if $\beta_j<\alpha_{j+1}$. 
In the second case $[\beta_j,\alpha_{j+1}]$ is called a non degenerate spectral gap. 
Notice also that, in some situations, the band $\Sigma_j$ 
can degenerate into a single point, i.e. $\alpha_{j}=\beta_j$. 
This happens if, and only if,  $\lambda_j(p)$   is  constant on $[-\pi,\pi]$.
In which case, this constant value is an infinitely degenerate eigenvalue of $H_0$. 
To illustrate these considerations we will give in Section \ref{libre} some explicit examples 
and show  necessary and/or sufficient conditions
to a spectral band to be non degenerate. 

Using  the eigenvalues $\lambda_j(p)$'s of the reduced operators $h(p)$, we will construct a conjugate operator
  for $H_0$.  
That is to say a self-adjoint operator 
$\mathbb{A}$ such that, for each compact real interval $J$ included in
$\Bbb{R}\setminus \kappa (H_0)$, we have
\begin{equation}\label{EM1}
E_{H_0}(J)[H_0,i\mathbb{A}]E_{H_0}(J)\geq c\;E _{H_0}(J)
\end{equation}
 for some $c>0$. For the definition of the commutator $[H_0,i\mathbb{A}]$ see Section \ref{rappel}.
In particular, we obtain  several important resolvent estimates for $H_0$, see Theorem \ref{rbvcontinuity}. 
We used a similar strategy in \cite{BoS} to construct conjugate operators from dispersion curves for 
perturbations of fibered systems in cylinders.
 
Next we prove that if  $A_n-A$ and $B_n-B$
  tend  to zero fast enough as $n$ tends to infinity, then  the Mourre estimate (\ref{EM1}) is preserved for $H$ up to a compact operator. 
 More precisely, we will show that, if (\ref{HP1}) below holds, then
  for the same  $\mathbb{A}, J$ and $c>0$  as in (\ref{EM1}),  
\begin{equation}\label{EM2}
E_{H}(J)[H,i\mathbb{A}]E_{H}(J)\geq c\;E _{H}(J)+K
\end{equation}
for some compact operator $K$. In particular, we get
\begin{theorem}\label{theorem0}
 Suppose that 
\begin{equation}\label{HP1}
\lim_{|n|\rightarrow\infty}|n|(\|A_n-A\|+\|B_n-B\|)=0.
\end{equation}
 Then outside $\kappa(H_0)$  the eigenvalues of $H$ are all finitely degenerate and their possible accumulation points are contained  in $\kappa(H_0)$. 
 \end{theorem}
Put $\langle {x}\rangle=\sqrt{1+x^2}$
for a real number $x$. We denote by $\textbf{N}$ the multiplication operator defined by 
$(\textbf{N}\psi)_n=n\psi_n$ for all $\psi\in\mathcal{H}$, 
and by   $B(\mathcal{H})$ the Banach algebra of bounded operators in the Hilbert space $\mathcal{H}$.
We also need the interpolation space $\mathcal{K}:=(D(\textbf{N}),\mathcal{H})_{1/2,1}$ which can be described, according to  Theorem 3.6.2 of  \cite{ABG}, by the norm, 
$$
\|\psi\|_{1/2,1}=\sum_{j=0}^\infty2^{j/2}\|\theta(2^{-j}|\textbf{N}|)\psi\|
$$ 
where $\theta\in C^\infty_0(\mathbb{R})$ with  $\theta(x)>0$ if $2^{-1}<x<2$ and $\theta(x)=0$ otherwise. Clearly,  
$\langle {\textbf{N}}\rangle^{-s}\mathcal{H}\subset\mathcal{K}\subset\mathcal{H}$ for any $s>1/2$.
\begin{theorem}\label{theorem1}
 Suppose that the perturbation $V$ satisfies
\begin{equation}\label{HP}
\int_1^\infty\|\theta(|\textbf{N}|/r)V\|dr<\infty.
\end{equation}
Then the limits
$
 ({H}-x\mp i0)^{-1}:=\lim_{\mu\rightarrow0+}({H}-x\mp i\mu)^{-1}
$
exist locally uniformly on  $\mathbb{R}\setminus[\kappa(H_0)\cup\sigma_p(H)]$ for the weak* topology of 
$B({\mathcal{K},\mathcal{K}^*})$.
In particular,  the singular continuous spectrum $\sigma_{sc}(H)$ is empty.
\end{theorem}
Mention that one may show that condition (\ref{HP}) is equivalent to
\begin{equation}\label{HP-explicit}
\int_1^\infty\sup_{r<|n|<2r}(\|A_n-A\|+\|B_n-B\|)dr<\infty.
\end{equation}\
Hence these conditions are clearly satisfied if, for some $\theta>0$, one has
\begin{equation}\label{HP'}
\sup_{n\in\mathbb{Z}}|n|^{1+\theta}(\|A_n-A\|+\|B_n-B\|)<\infty.
\end{equation}
Remark also that  Theroem \ref{theorem1} implies that 
$
 \langle {\textbf{N}}\rangle^{-s}({H}-x\mp i0)^{-1}\langle {\textbf{N}}\rangle^{-s}\in B(\mathcal{H})
 $, 
for any $s> \frac{1}{2}$.
In our next theorem we describe continuity properties of these boundary values of the resolvent of $H$ 
as functions of $x$, as well as some of their propagation consequences. For we need the following definition.
Let $(E,\|\cdot\|)$ be a Banach space and $f:{\Bbb R}\rightarrow E$ be a bounded
continuous function.  For an integer $m\geq1$ let $w_m$ be the modulus of
continuity of order $m$ of $f$ defined on $(0,1)$ by
$$
w_m(f,\varepsilon)=\sup_{x\in\mathbb{R}}
\biggl\|\sum_{j=0}^m (-1)^j
\left(
\begin{array}{cc}
 m  \\ j   
\end{array}
\right)
f(x+j\varepsilon)\biggr\|.
$$
We say that $f\in\Lambda^{\alpha}$, $\alpha>0$,
 if  there is an integer $m>\alpha$ such
that 
$$
\sup_{0<\varepsilon<1}\varepsilon^{-\alpha} w_m(\varepsilon)<\infty.
$$
Notice that if $\alpha\in(0,1)$ then $\Lambda^\alpha$ is nothing but the space of H\H{o}lder continuous functions of order $\alpha$.
In contrast, if $\alpha=1$ then $\Lambda^1$ consists of smooth functions in Zygmund's sense (they are not Lipschitz in general, see \cite{Z}). Finally, if $\alpha>1$ and 
$n_\alpha$ is the greatest integer strictly less than $\alpha$  then $f$ belongs to $\Lambda^\alpha$ if, and only if, $f$ is $n_\alpha-$times continuously differentiable with bounded derivatives and its derivative
 $f^{(n_\alpha)}$ of order $n_\alpha$  is of class $\Lambda^{\alpha-n_\alpha}$. 
 For example, $f$ belongs to $\Lambda^2$  means that $f$ is continuously differentiable with a bounded derivative and $f'$ is of class  $\Lambda^1$ ($f'$ is not Lipschitz in general).
%%%%%%%%%%%%%%%%%%%%%%%%%%%%%%%%%%%%%%%%%%%%%%%%%%%%%%%%%%%%%%%%%%
%%%%%%%%%%%%%%%%%%%%%%%%%%%%%%  Theorem 1.2  %%%%%%%%%%%%%%%%%%%%%
%%%%%%%%%%%%%%%%%%%%%%%%%%%%%%%%%%%%%%%%%%%%%%%%%%%%%%%%%%%%%%%%%%
\begin{theorem}\label{theorem11}
Suppose that (\ref{HP'}) holds. Then, for any $ \frac{1}{2}<s\leq \frac{1}{2}+\theta$,  the maps
$$
 x\longmapsto  \langle {\textbf{N}}\rangle^{-s}({H}-x\mp i0)^{-1}\langle {\textbf{N}}\rangle^{-s}\in B({\mathcal   H})
$$
are locally of class  $\Lambda^{s- \frac{1}{2}}$ on $\mathbb{R}\setminus[\kappa(H_0)\cup\sigma_p(H)]$. 
Moreover, 
$$
||\langle {\textbf{N}}\rangle^{-s}e^{-i{H}t}\varphi(H)\langle {\textbf{N}}\rangle^{-s}||\leq C(1+|t|)^{-(s- \frac{1}{2})}.
$$
for all $\varphi\in C^\infty_0(\mathbb{R}\setminus[\kappa(H_0)\cup\sigma_p(H)])$ and some $C>0$.
\end{theorem}
%%%%%%%%%%%%%%%%%%%%%%%%%%%%%%%%%%%%%%%%%%%%%%%%%%%%%%%%%%%%%%%%%%
%%%%%%%%%%%%%%%%%%%%%%%%%%%%%%  Theorem 1.1  %%%%%%%%%%%%%%%%%%%%%
Notice that for $H_0$ we have more than that, see Theorem \ref{rbvcontinuity}. Moreover, Theorems \ref{theorem0}-\ref{theorem11}  are actually valid for operators of the form $H=H_0+V$ where $V$ is a symmetric compact perturbation such that
 $\langle {\textbf{N}}\rangle^{1+\theta}V$ is bounded with $\theta>0$.  See section \ref{perturbŽ} where this remark is discussed more precisely.

%%%%%%%%%%%%%%%%%%%%%%%%%%%%%%%%%%%%%%%%%%%%%%%%%%%%%%%%%%%%%%%%%%
Finally, in the second part of the paper, we apply these abstract results to some concrete models. 
The first application concerns the case where   $A$ is positive definite.
It is one of the most common assumptions in the studies
of block Jacobi matrices and their applications to matrix orthogonal polynomial theory.
In this case,  all the eigenvalues of $h(p)$ are monotonic functions of $p$. 
 This applies directly to a	 class of difference operators on cylindrical domains such as discrete wave propagators. 
 We devote the next applications to cases where the matrix $A$ is singular. 
Then we show how  Jacobi operators with periodic coefficients fit in this situation. 
Notice that while the spectral results of Jacobi operators with periodic coefficients are well known
our Mourre estimate is a new ingredient. Finally, we study a model where $A$ is a lower triangular matrix that 
we apply  to  symmetric difference operators of fourth-order.

The paper is organized as follows.  
\begin{itemize}
\item In Section \ref{libre} we study in details the unperturbed operator $H_0$. 
\item Section \ref{rappel} contains a brief review on what we need from Mourre's theory.  
\item In Sections \ref{multiplication} and \ref{mestimateH0} we construct a conjugate operator for $H_0$. 
\item In section \ref{perturbŽ} we  extend our Mourre estimate to the operator $H$. 
\item In section \ref{positivedefinie} we study the case where $A$ is positive definite.
\item Section  \ref{waves} is devoted  to  discrete wave propagators in stratified cylinders. 
\item  In  Sections \ref{app-singular1}-\ref{singular3} we discuss cases where $A$ is singular.
\item In section \ref{triangulaire} we consider
 a case where $A$ is a lower triangular.
 \item Finally  in Section \ref{further} we explain how our method extends to difference operator of higher order with matrix coefficients. 
\end{itemize}
 
 \textbf{Acknowledgements:  }{ We would like to thank   J. Janas and  P. Cojuhari
for their hospitality and useful discussions on the subject during 
my visit to IPAN (Krakow).}
\protect\setcounter{equation}{0}
%
%
%%%%%%%%%%%%%%%%%%%%%%%%%%%%%%%%%%%%%%%%%%%%%%%%%%%%%%%%%%%%%%%%%%
%%%%%%%%%%%%%%%%%%%%%%%%%%%%%%%%%%%%%%%%%%%%%%%%%%%%%%%%%%%%%%%%%%
%%%%%%%%%%%%%%%%%%%%%%%%   Section 3 %%%%%%%%%%%%%%%%%%%%%%%%%%%%%
%%%%%%%%%%%%%%%%%%%%%%%%%%%%%%%%%%%%%%%%%%%%%%%%%%%%%%%%%%%%%%%%
%%%%%%%%%%%%%%%%%%%%%%%%%%%%%%%%%%%%%%%%%%%%%%%%%%%%%%%%%%%%%%%%
%%%%%%%%%%%%%%%%%%%   L'op{\'e}rateur non perturb{\'e}    %%%%%%%%%%%%%%
%%%%%%%%%%%%%%%
%%%%%%%%%%%%%%%%%%%%%%%%%%%%%%%%%%%%%%%%%%%%%%%%%
%%%%%%%%%%%%%%%%%%%%%%%%%%%%%%%%%%%%%%%%%%%%%%%%%%%%%%%%%%%%%%%%
\setcounter{equation}{0}

%%%%%%%%%%%%%%%%%%%%%%%%%%%%%%%%%%%%%%%%%%%%%%%%%%%%%%%%%%%%%%%%%%%%%%%%%
%%%%%%%%%%%%%%%%%%%%%%%%%%%%%%%%%%%%%%%%%%%%%%%%%%%%%%%%%%%%%%%%%%%%%%%%%
\protect\setcounter{equation}{0}
%%%%%%%%%%%%%%%%%%%%%%%%%%%%%%%%%%%%%%%%%%%%%%%%%%%%%%%%%%%%%%%%
%%%%%%%%%%%%%%%%%%%%%%%%%%%%%%%%%%%%%%%%%%%%%%%%%%%%%%%%%%%%%%%%
\section{Basic properties of $H_0$}\label{libre}
%%%%%%%%%%%%%%%%%%%%%%%%%%%%%%%%%%%%%%%%%%%%%%%%%%%%%%%%%%%%%%%%
%%%%%%%%%%%%%%%%%%%%%%%%%%%%%%%%%%%%%%%%%%%%%%%%%%%%%%%%%%%%%%%%
%%%%%%%%%%%%%%%%%%%%%%%%%%%%%%%%%%%%%%%%%%%%%%%%%%%%%%%%%%%%%%%%
%%%%%%%%%%%%%%%%%%%%%%%%%%%%%%%%%%%%%%%%%%%%%%%%%%%%%%%%%%%%%%%%
Consider the unitary transform 
$\mathcal{F}:l^2(\mathbb{Z},\mathbb{C}^N)\rightarrow L^2([-\pi,\pi],\mathbb{C}^N)$ defined by
$$
\mathcal{F}(\psi)(p)=\widehat\psi(p)=\sum_{n\in\mathbb{Z}}e^{inp}\psi_n~~,\mbox{ for each }~~\psi\in l^2(\mathbb{Z},\mathbb{C}^N).
$$
Recall that the unperturbed operator $H_0$ is defined in $\mathcal{H}$   by
$$
(H_0\psi)_n=A^*\psi_{n-1}+B\psi_n+A\psi_{n+1}, ~~\mbox{
for all } n\in\mathbb{Z}.
$$
Direct computation shows that,  for each $\psi\in l^2(\mathbb{Z},\mathbb{C}^N)$,
$$
\widehat{H_0\psi}(p)=(e^{-ip}A+e^{ip}A^*+B)\widehat{\psi}(p)=h(p)\widehat{\psi}(p)~~,\mbox{ for all }~p\in[-\pi,\pi]~.
$$
Hence (\ref{direct})-(\ref{reduced}) are proved. 
Recall that, for all $p\in\Bbb{R}$,  $\{\lambda_{j}(p);\,1\leq j\leq N\}$ are the repeated eigenvalues  of
$h(p)$  and $\{W_{j}(p);\,1\leq j\leq N\}$ is a corresponding orthonormal basis
of eigenvectors.
The Proposition \ref{prop1} follows from:
\begin{proposition}\label{prop3.1}
The operator $H_0$ is unitarily equivalent to $M=\oplus_{j=1}^{j=N}\lambda_j(p)$ acting in 
the direct sum $\oplus_{j=1}^{j=N}L^2([-\pi,\pi])$. 
\end{proposition}
\textbf{Proof } Consider the operator 
$
\mathcal{U}:l^2(\mathbb{Z},\mathbb{C}^N)\rightarrow\oplus_{j=1}^{j=N}L^2([-\pi,\pi])
$
defined by 
$
\mathcal{U}\psi=(f_j)_{1\leq j\leq N}
$,
where 
\begin{equation}\label{deftrsfF}
f_j(p)=\langle W_j(p),\widehat\psi(p)\rangle_{\mathbb{C}^N}=\sum_{n\in\mathbb{Z}}\langle W_j(p),\psi_n\rangle_{\mathbb{C}^N}e^{inp}.
\end{equation}
Clearly, $\mathcal{U}$ 
is a unitary operator with the following inversion formula:
$$
(\mathcal{U}^{-1}f)_n=\sum_{j=1}^N\frac{1}{2\pi}\int_{-\pi}^\pi f_j(p)W_j(p)e^{-inp}dp,
$$
for all $f=(f_1,\cdots,f_N)$ in $\oplus_{j=1}^{j=N}L^2([-\pi,\pi])$.
Moreover, direct computation shows that,  for each $\psi\in l^2(\mathbb{Z},\mathbb{C}^N)$,
$$
(\mathcal{U}H_0\psi)_j(p)=\langle W_j(p),\widehat{H_0\psi}(p)\rangle_{\mathbb{C}^N}=\langle W_j(p),h(p)\widehat{\psi}(p)\rangle_{\mathbb{C}^N}=\lambda_j(p)f_j(p).
$$
The proof is complete.

In the sequel  we denote the imaginary part of a complex number 
 $z$ by $\Im z$.
%%%%%%%%%%%%%%%%%%%%%%%%%%%%%%%%%%%
%%%%%%%%%%%%%%%%%%%%%%%%%%%%%%%%%%%
\begin{proposition}\label{deŽrivŽee}
The eigenvalue $\lambda_j(p)$ is identically constant on $[-\pi,\pi]$ if, and only if, 
$$
\lambda_j'(p)=\Im(e^{ip}\langle AW_j(p),W_j(p)\rangle_{\mathbb{C}^N})=0, ~~\mbox{ for all  }~
p\in[-\pi,\pi].
$$
 In this case, the spectral band $\Sigma_j=\lambda_j([-\pi,\pi])$ degenerates into a single point 
 which is an infinitely degenerate eigenvalue of $H_0$.
\end{proposition}
%%%%%%%%%%%%%%%%%%%%%%%%%%%%%%%%%%%
%%%%%%%%%%%%%%%%%%%%%%%%%%%%%%%%%%%
\textbf{Proof:  } We know that,  for all  $p\in[-\pi,\pi]$,
$$
\lambda_j(p)=\langle h(p)W_j(p),W_j(p)\rangle_{\mathbb{C}^N}.
$$  
Since $\|W_j(p)\|^2=1$ then 
$
\langle W_j'(p),W_j(p)\rangle_{\mathbb{C}^N}+\langle W_j(p),W'_j(p)\rangle_{\mathbb{C}^N}=0.
$
Thus
\begin{eqnarray*}
\lambda_j'(p)=\langle h(p)'W_j(p),W_j(p)\rangle_{\mathbb{C}^N}
&=&\langle i(-e^{-ip}A+e^{ip}A^*)W_j(p),W_j(p)\rangle_{\mathbb{C}^N}\\
&=&-2\Im\big(e^{ip}\;\langle AW_j(p),W_j(p)\rangle_{\mathbb{C}^N}\big).
\end{eqnarray*}
The proof is  complete.

\textbf{Remark } In particular, the j$^{th}$ spectral band $\Sigma_j=\lambda_j([-\pi,\pi])$ is not degenerate if, and only
 if, there is $p\in[-\pi,\pi]$ such that
$
\lambda_j'(p)=-2 \Im\big(e^{ip}\;\langle AW_j(p),W_j(p)\rangle_{\mathbb{C}^N}\big)\not=0
$. 
In this case,
  $\sigma_p(H_0)\cap \Sigma_j$ is finite and $\sigma_{sc}(H_0)\cap \Sigma_j=\emptyset$. 
 Moreover, $\Sigma_j$ contains an eigenvalues $\mathbf{e}$ 
 of $H_0$ if and only if there exists $i\not=j$ such that $\lambda_{i}(q)=\mathbf{e}$ for all $q\in[-\pi,\pi]$.
Finally, if for all $ j=1,\cdots,N$ there is $p\in[-\pi,\pi]$ such that $\lambda'_j(p)\not=0$,
then the spectrum of $H_0$ is purely absolutely continuous, i.e. $\sigma_p(H_0)=\sigma_{sc}(H_0)=\emptyset$. We 
deduce immediately the following corollary:
%
%%%%%%%
%%%%%%%Le cas ou $A$ est définie positive
%%%%%%%
 \begin{corollary}\label{cas-positive}
If $A$ is positive definite then  $\lambda_1(p),\cdots,\lambda_N(p)$  are even on $[-\pi,\pi]$, decreasing on $[0,\pi]$
and  for  all  $ j=1,\cdots,N$, we have
 $$
\lambda_j'(p)=-2\langle AW_j(p),W_j(p)\rangle_{\mathbb{C}^N}\sin p~,~~\mbox{  for all }~p\in[-\pi,\pi].
$$
In particular, $\Sigma_j=[\lambda_j(\pi),\lambda_j(0)]$, $\kappa(\lambda_j)=\{\lambda_j(\pi),\lambda_j(0)\}$ and  
 the spectrum of $H_0$ is purely absolutely continuous with
$
\sigma(H_0)=\cup_{1\leq j\leq N}[\lambda_j(\pi),\lambda_j(0)].
$
\end{corollary} 
%
%
%%%%%%%
%%%%%%% Le cas ou $A$ est inversible symétrique
%%%%%%%
 \begin{corollary}\label{non-singular}
If $A$ is a symmetric invertible matrix then  
 $\lambda_1(p), \cdots ,\lambda_N(p)$  are non constant even functions on $[-\pi,\pi]$.
  In particular, the spectrum of $H_0$ is purely absolutely continuous.
  \end{corollary} 
\textbf{Proof } Assume by contradiction that $\lambda_1(p)=c$ for all  $p\in [-\pi,\pi]$ for some constant $c$.
Then for any $p\in [-\pi,\pi]$ there exists $W(p)$ an eigenvector of $h(p)$ associated to $c$ so that $h(p)W(p)=cW(p)$.
Hence for any $p\in [-\pi,\pi]$, $(2\cos{p})AW(p)=[(c-B)A^{-1}]AW(p)$,  which means that  $[-2,2]\subset\sigma((c-B)A^{-1})$ which is impossible.
%
%
%
%%%%%%%
%%%%%%% Le cas ou $A$ est  symétrique singulière
%%%%%%%
 \begin{corollary}\label{singulière}
 Assume that $A$ is a symmetric matrix such that for any eigenvector $w$ of $B$,
 $\langle A w,w\rangle\not=0$. Then the conclusion of  Corollary \ref{non-singular} holds true.
 \end{corollary} 
\textbf{Proof }  Here remark that $h(\pi/2)=B$. So $\lambda_j(\pi/2)$ is an eigenvalue of $B$ and $W_j(\pi/2)$ is an associated eigenvector.
So $\lambda'_j(\pi/2)=-2\langle A W_j(\pi/2),W_j(\pi/2)\rangle\not=0$, by hypothesis. The proof is complete.

%
%
%
%%%%%%%
%%%%%%% Le cas ou $A$ est  non symétrique 
%%%%%%%
 \begin{corollary}\label{nonsymmetric}
 Assume that  for any eigenvector $w$ of $\;i(A-A^*)+B$ one has
 $\langle (A+A^*) w,w\rangle\not=0$. Then  
  $\lambda_1(p),\cdots,\lambda_N(p)$  are non constant  on $[-\pi,\pi]$. In particular, 
 the spectrum of $H_0$ is purely absolutely continuous and there exist 
 $\alpha_1<\beta_1, \cdots,\alpha_n<\beta_n$
such that  
$
\sigma(H_0)=\cup_{1\leq j\leq N}[\alpha_j,\beta_j].
$
\end{corollary} 
\textbf{Proof } It is enough to remark that $h(-\pi/2)=i(A-A^*)+B$.

We close this section by few  explicit examples to illustrate  these consideration and
our comment just after Proposition \ref{prop1} on the spectrum of $H_0$.
\begin{example} 
Assume that $A$ and $B$  are both diagonals with diagonal elements  $A_{jj}=a_j\geq0 $ and $B_{jj}=b_j\in\mathbb{R}$.
In this case, the eigenvalues of the reduced operators $h(p)$ are given by $\lambda_j(p)=b_j+2a_j\cos p$. 
Therefore, we immediately see that by playing with the parameters $a_j$ and $b_j$ one may realize whatever we said 
with the associated spectral bands $\Sigma_j=[b_j-2a_j,b_j+2a_j]$.
\end{example}
\begin{example} 
Let 
$
A=\left(
\begin{array}{cc}
  0&1     \\
  1&  0   
\end{array}
\right)\quad \mbox{and}\quad B=\left(
\begin{array}{cc}
  b&a   \\
  a&  -b
\end{array}
\right)$
where $a,b>0$.
 In  this case, the eigenvalues of  $h(p)$  are given, for all $p\in[-\pi,\pi]$, by
$$
\lambda_j(p)=(-1)^j\sqrt{b^2+(a +2\cos p)^2}\quad \mbox{with } j=1,2.
$$
Put $\alpha=\sqrt{b^2+(a -2)^2}$ and $\beta=\sqrt{b^2+(a +2)^2}$. Hence we have the following.
\begin{enumerate}
\item If $a\geq2$ then $\lambda_1$ is increasing on $[0,\pi]$ while $\lambda_2$ is decreasing and
$
\kappa(H_0)=\{\lambda_j(0), \lambda_j(\pi)\}_{j=1,2}
$.  Moreover, 
 the spectrum of $H_0$ is purely absolutely continuous and consists of  two spectral bands with a  gap between them of length $2\alpha$ and 
 $
\sigma(H_0)=[-\beta,-\alpha]\cup[\alpha,\beta]
$. 
\item If $0<a<2$ then the critical points of $\lambda_1$ and $\lambda_2$ are $0,\pi$ plus an additional  
point $p_0\in]\pi/2,\pi[$ defined by $\cos(p_0)=-a/2$ . In such case, 
$
\sigma(H_0)=[-\beta,-b]\cup[b,\beta]
$
and 
$
\kappa(H_0)=\{\lambda_j(0), \lambda_j(p_0),\lambda_j(\pi)\}_{j=1,2}
$.
So we have  two spectral bands with a gap between them of length $2b$. 
\end{enumerate}
\end{example}
\begin{example}  
Let 
$
A=\left(
\begin{array}{cc}
  0&0     \\
  a_2&  0   
\end{array}
\right) \mbox{,} B=\left(
\begin{array}{cc}
  b&a_1   \\
  a_1& -b
\end{array}
\right)$
where $a_1,a_2>0$ and $b\in\mathbb{R}$. 
 In  this case, the eigenvalues of reduced operators $h(p)$   are   given by
$$
\lambda_j(p)=(-1)^j\sqrt{b^2+a_1^2+a_2^2+2a_1a_2\cos p}~,\quad j=1,2.
$$
 Then  $\lambda_1(p)$ is an increasing function of $p\in[0,\pi]$ while $\lambda_2(p)$ is decreasing and $\lambda_1(p)<\lambda_2(p)$, for all $p\in]0,\pi[$. Therefore,
$\kappa(H_0)= \{\lambda_1(0),\lambda_1(\pi),\lambda_2(\pi),\lambda_2(0)\}$ and 
$
\sigma(H_0)=[\lambda_1(0),\lambda_1(\pi)]\cup[\lambda_2(\pi),\lambda_2(0)].
$
This spectrum is purely absolutely continuous and consists of  two spectral bands with a gap between them of length
$
L=\lambda_2(\pi)-\lambda_1(\pi)=2\sqrt{b^2+(a_1-a_2)^2}.
$
That gap is degenerate (i.e. $L=0$), if and only if, $b=0$ and $a_1=a_2$.
In such case, the two spectral bands have one common point 
 $\lambda_1(\pi)=\lambda_2(\pi)=2a$ which is a double eigenvalue of $h(\pi)$.

Mention that $\lambda_1$ and $\lambda_2$ are identically constant on $[-\pi,\pi]$  if, and only if,   $a_1a_2=0$. 
For example, if $a_2>0$ and $a_1=0$
then the spectrum of $H_0$ consists of two infinitely degenerate eigenvalues  $\pm\sqrt{b^2+a_2^2}$.
\end{example}
%%%%%%%%%%%%%%%%%%%%%%%%%%%%%%%%%%%%%%%%%%%%%%%%%%%%%%%%%%%%%%%%%%
%%%%%%%%%%%%%%%%%%%%%%%%%%%%%%%%%%%%%%%%%%%%%%%%%%%%%%%%%%%%%%%%%%
%%%%%%%%%%%%%%%%%%%%%%%%%%%%%%%%%%%%%%%%%%%%%%%%%%%%%%%%%%%
%%%%%%%%%%%%%%%%%%%%%%%%%%%%%%%%%%%%%%%%%%%%%%%%%%%%%%%%%%%
%%%%%%%%%%%%%%%%%%%%%%%%%%%%%%%%%%%%%%%%%%%%%%%%%%%%%%%%%%%%%%%%%%%
%            THE CONJUGATE OPERATOR METHOD
%%%%%%%%%%%%%%%%%%%%%%%%%%%%%%%%%%%%%%%%%%%%%%%%%%%%%%%%%%%%%%%%%%%
%
\section{ The conjugate operator theory}\label{rappel}
The following brief review on the conjugate operator theory is based on
\cite{ABG,BG,BGS2,S1}. 
Let $\mathbb{A}$ be a  self-adjoint operator in a  separable complex Hilbert space  $\mathcal{H}$ and $S\in B({\mathcal H})$.
%%%%%%%%%%%%%%%%%%%%%%%%%%%%%%%%%%%%
\begin{definition} \label{defreg}
(i) Let $k\geq 1$ be an integer and $\sigma>0$.
 We say that  $S$ is of class $C^{k}(\mathbb{A})$, respectively of class $\mathcal{C}^{\sigma}(\mathbb{A})$, if the map
$$
t\longmapsto S(t)= e^{-i\mathbb{A}t}Se^{i\mathbb{A}t}\in B({\mathcal H})
$$
is strongly of class $C^k$, respectively of class $\Lambda^\sigma$ on ${{\Bbb R}}$.

(ii)  We say that $S$ is of class $\mathcal{C}^{1,1}(\mathbb{A})$ if
$$
\int_0^1\|e^{-i\mathbb{A}\varepsilon}Se^{i\mathbb{A}\varepsilon}-2S+e^{i\mathbb{A}\varepsilon}Se^{-i\mathbb{A}\varepsilon}S\|\frac{d\varepsilon}{\varepsilon^2}<\infty.
$$
\end{definition}
\textbf{Remarks} 
(i) We have the following inclusions 
$\mathcal{C}^s(\mathbb{A})\subset \mathcal{C}^{1,1}(\mathbb{A})\subset C^1(\mathbb{A})$, for any $s>1$.

(ii) 
One may show  that $S$ is  of class $C^1(\mathbb{A})$, if and only if, the sequilinear form defined on 
$D(\mathbb{A})$ by 
$[S,\mathbb{A}]=S\mathbb{A}-\mathbb{A}S$  has a continuous extension to $\mathcal H$, which we identify with
the associated bounded operator in $\mathcal H$ (from the Riesz Lemma) that
we denote  by the same symbol.  Moreover,
$
[S,i\mathbb{A}]=\frac{d}{dt}|_{t=0}S(t).
$

(iii) Recall that  $\langle\mathbb{A}\rangle=(1+\mathbb{A}^2)^{1/2}$. To prove that $S$  is of class 
$\mathcal{C}^{\sigma}(\mathbb{A})$ it is enough to show that
  $\langle\mathbb{A}\rangle^\sigma S$ is bounded in $\mathcal H$,
see for example the appendix of \cite{BS1}.  In particular,  in the case where $\sigma=1+\theta$ for some $\theta>0$, the operator $S$ is of class ${C}^{1+\theta}(\mathbb{A})$ if one of the following conditions is true:\\
(1) $\langle\mathbb{A}\rangle^\sigma S$ is bounded in $\mathcal H$ or\\
(2) $S$ is of class  ${C}^{1}(\mathbb{A})$ and $\langle\mathbb{A}\rangle^\theta [S,i\mathbb{A}]$ is bounded in $\mathcal H$.

In the sequel of this section, let   $H$ be a bounded self-adjoint operator which is at least of class $C^1(\mathbb{A})$. 
Then $[H,i\mathbb{A}]$ defines  a   bounded operator
 in $\mathcal{H}$ that we still denote by the
 same symbol $[H,i\mathbb{A}]$.

We  define the open set $\tilde\mu^{\mathbb{A}}(H)$ of real numbers $x$ such that,
for some constant $a>0$, a neighborhood $\Delta$ of $x$ and  a
compact operator $K$ in $\mathcal H$, we have
\bea\label{mourre estimate}
E_H(\Delta )[H,i\mathbb{A}]E_H(\Delta )\geq aE_H(\Delta )+K.
\eea
The inequality (\ref{mourre estimate}) is called {\it the Mourre estimate} and the set of point 
$x\in\tilde\mu^{\mathbb{A}}(H)$ where it holds with $K=0$ will be denoted by $\mu^{\mathbb{A}}(H)$.
\begin{theorem} \label{virial} The set
$\tilde\mu^{\mathbb{A}}(H)$ contains at most a discrete set of eigenvalues
of $H$ and all these eigenvalues are finitely degenerate. Moreover,
$
\mu^{\mathbb{A}}(H)=\tilde\mu^{\mathbb{A}}(H)\setminus\sigma_{p}(H).
$
\end{theorem}
Recall that $\mathcal{K}_\mathbb{A}:=\mathcal{H}_{1/2,1}$ is the Besov space associated to 
$\mathbb{A}$ defined by the norm
$$
\|\psi\|_{1/2,1}=\|\tilde\theta(\mathbb{A}) \psi\|+\sum_{j=0}^\infty2^{j/2}\|\theta(2^{-j}|\mathbb{A}|)\psi\|.
$$
where $\theta\in C^\infty_0(\mathbb{R})$ with  $\theta(x)>0$ if $2^{-1}<x<2$ and $\theta(x)=0$ otherwise, and $\tilde\theta\in C^\infty_0(\mathbb{R})$ with  $\tilde\theta(x)>0$ if $|x|<2$ and $\tilde\theta(x)=0$ otherwise.
One has, see \cite{BG,S1}:
\begin{theorem} \label{lap}
 Assume that $H$ is of class $\mathcal{C}^{1,1}(\mathbb{A})$. Then $H$ has no singular continuous spectrum in $\mu^{\mathbb{A}}(H)$. 
 Moreover, the  limits
$
({H}-x\mp i0)^{-1} :=
\lim_{\mu\rightarrow0+}({H}-x\mp i\mu)^{-1}
$
exist locally uniformly on  $\mu^{\mathbb{A}}(H)$ for the weak* topology of $B(\mathcal{K}_\mathbb{A},\mathcal{K}_\mathbb{A}^*)$.
\end{theorem}
We also have the following:
\begin{theorem} \label{lap1}
 Assume that $H$ is of class $\mathcal{C}^{s+\frac{1}{2}}(\mathbb{A})$
for some $s>1/2$. Then the  maps
$
 x\mapsto\langle\mathbb{A}\rangle^{-s}({H}-x\mp i0)^{-1}\langle\mathbb{A}\rangle^{-s}
$
are locally of class
$\Lambda^{s-\frac{1}{2}}$ on $\mu^{\mathbb{A}}(H)$.
Moreover, for every $\varphi\in C^\infty_0(\mu^{\mathbb{A}}(H))$ we have
$$
||\langle\mathbb{A}\rangle^{-s}e^{-iHt}\varphi(H)\langle\mathbb{A}\rangle^{-s}||\leq C<t>^{-(s-\frac{1}{2})}.
$$
\end{theorem}
Since $\|e^{-iHt}\varphi(H)\|=\|\varphi(H)\|$, then by using  the complex interpolation one  obtains much more than this, see \cite{BG,BGS2,S1}.
For example, if $H$ is of class $C^\infty(\mathbb{A})$, then, for any  $\varphi\in C^\infty_0(\mu^{\mathbb{A}}(H))$ and  $\sigma,\varepsilon>0$, we have
$$
||\langle\mathbb{A}\rangle^{-\sigma}e^{-i{H}t}\varphi(H)\langle\mathbb{A}\rangle^{-\sigma}||\leq C(1+|t|)^{-\sigma+\varepsilon},\quad \mbox{for all }t\in\mathbb{R}
$$
\protect\setcounter{equation}{0}
%

%%%%%%%%%%%%%%%%%%%%%%%%%%%%%%%%%%%
\section{Mourre estimate for $\lambda _{j}$}\label{multiplication}
%%%%%%%%%%%%%%%%%%%%%%%%%%%%%%%%%%%

Here we denote  the operator of multiplication by a function 
$f$ in $L^2([-\pi,\pi])$  by the same symbol  $f$ or 
by $f(p)$ when we want to stress the $p$-dependence. 

According to Section \ref{libre},  $H_0$ is unitarily equivalent to $\oplus_{j=1}^N\lambda_j$. 
So it is convenient to prove first the  Mourre estimate for each multiplication operator $\lambda_j$.
So, let us  fix $j\in\{1,\cdots,N\}$ and put $\lambda=\lambda_j$. 
Recall that the set  $\kappa (\lambda)$ of critical values of $\lambda$ is finite. 
Let $\Delta\subset \lambda([-\pi,\pi])$ be
a compact set, such that ${\Delta}\cap \kappa (\lambda)=\emptyset.$
Obviously, on the set $\lambda^{-1}(\Delta)$ one has $|\lambda ^{\prime }(p)|^2\geq c>0$.
One may then construct a function $F\in C^{\infty } ({[-\pi,\pi]}),$
such that 
\begin{equation}\label{FF}
F(p) \lambda^{\prime }(p)\geq c>0\quad\mbox{ on $\lambda^{-1}(\Delta)$.}
\end{equation}  
For example,  $F=\lambda'$   perfectly fulfills this role. Nevertheless, 
in some concrete situations  more appropriate choices can be made, see next sections.
We are now able to define our conjugate operator $\mathbf{a}$ to $\lambda$ on $L^2([-\pi,\pi])$ by,
 %%%%%%%%%%%%%%%%
\begin{equation}
\mathbf{a}=  \frac{i}{2}\{F(p) \frac{d}{dp}+ \frac{d}{dp} F(p)\}=
i F(p) \frac{d}{dp} +\frac{i}{2}F^{\prime }(p). \label{ocdef}
\end{equation}
%%%%%%%%%%%%%%%%
It is clear that $\mathbf{a}$ is an essentially self-adjoint operator  in $L^2([-\pi,\pi])$. 
Moreover, direct computation yields to
$$
[\lambda,i\mathbf{a}]=F(p) \lambda^{\prime }(p).
$$
So $\lambda$ is of class $C^{1}(\mathbf{a})$,  and by repeating the same calculation, we show that  
$\lambda$ is of class $C^{\infty}(\mathbf{a})$. Moreover, thanks to (\ref{FF}), $\mathbf{a}$ is strictly conjugate
to $\lambda$ on $\Delta$:
\begin{eqnarray}\label{mel}
E_{\lambda}(\Delta)[\lambda,i\mathbf{a}]E_\lambda{(\Delta)}
&=&E_{\lambda}(\Delta)\big(F(p) \lambda^{\prime }(p)\big)E_{\lambda}(\Delta)
\geq cE_{\lambda}(\Delta)
\end{eqnarray}
From now on we put $F=\lambda'$. We have:
\begin{proposition}\label{prop3.3}
The operator $\lambda$ is of class $C^{\infty}(\mathbf{a})$ and $\mathbf{a}$ is locally  strictly conjugate
to $\lambda$ on $\mathbb{R}\setminus\kappa(\lambda)$, i.e.
$
\mu^{\mathbf{a}}(\lambda)=\mathbb{R}\setminus\kappa(\lambda).
$
\end{proposition}
\textbf{Proof  }
First, recall that $\mathbb{R}\setminus\sigma(\lambda)\subset \mu^{\mathbf{a}}(\lambda)$. So,
 if $\lambda$ is constant then the assertion is trivial, since $\sigma(\lambda)=\kappa(\lambda)$.  
 Suppose that  $\lambda$ is non constant and
let $x\in \sigma(\lambda)\setminus\kappa(\lambda)$.  Then choose a compact interval $\Delta$ about $x$
 sufficiently small so that $\Delta\cap \kappa (\lambda)=\emptyset.$
The last discussion  shows that $\mathbf{a}$ is strictly conjugate
to $\lambda$ on $\Delta$, i.e.  $x\in\mu^{\mathbf{a}}(\lambda)$. The proof is complete.

\protect\setcounter{equation}{0}
%

%%%%%%%%%%%%%%%%%%%%%%%%%%%%%%%%%%%%%%%
%%%%%%%%%%%%%%%%%%%%%%%%%%%%%%%%%%%%%%%
%%%%%%%%%%%%%%%%%%%%%%%%%%%%%%%%%%%%%%%
\section{Mourre estimate for $H_0$}\label{mestimateH0}
%%%%%%%%%%%%%%%%%%%%%%%%%%%%%%%%%%%%%%%
%%%%%%%%%%%%%%%%%%%%%%%%%%%%%%%%%%%%%%%
%%%%%%%%%%%%%%%%%%%%%%%%%%%%%%%%%%%%%%% 
 The direct sum $%
M=\oplus _{1\leq j\leq N}\lambda _{j}$ is a self-adjoint bounded operator in the Hilbert
space 
$%
\oplus _{1\leq j\leq N} L^{2}([-\pi,\pi])
$.
Recall that  $\kappa (H_0)= \cup
_{1\leq j\leq N}\kappa (\lambda _{j})$ is finite.
For any $j=1,\cdots,N$,  let us denote by $\mathbf{a}_{j}$ the operator given by (\ref{ocdef}),  with $F=F_{j}$ constructed in the last section with $\lambda=\lambda_j$,  that is
\begin{equation}
\mathbf{a}_j= \frac{i}{2}\{F_j(p) \frac{d}{dp}+ \frac{d}{dp} F_j(p)\}=
i F_j(p)\frac{d}{dp} +\frac{i}{2}F_j^{\prime }(p). \label{ocdef1}
\end{equation}
and define the operator 
\begin{equation}
\widetilde{A}=\oplus _{1\leq j\leq N}\mathbf{a}_{j} .  \label{defdeatld}
\end{equation}
\begin{proposition}\label{moureM} 
The operator $M$ is of class $C^{\infty}(\widetilde{A})$ and $\widetilde{A}$ is locally  strictly conjugate
to $M$ on  $\mathbb{R}\setminus\kappa(H_0)$, i.e.
$
\mu^{\widetilde{A}}(M)=\mathbb{R}\setminus\kappa(H_0).
$
\end{proposition}
\textbf{Proof:  }  Let $x\in\mathbb{R}\setminus\kappa(H_0)$. By definition $\kappa(H_0)=\cup_{1\leq j\leq N}\kappa(\lambda_j)$,
and according to  Proposition \ref{prop3.3}, $\mu^{{\mathbf{a}_j}}(\lambda_j)=\mathbb{R}\setminus\kappa(\lambda_j).$
Hence,  $x\in\cap_{1\leq j\leq N}\mu^{{\mathbf{a}_j}}(\lambda_j)$.  Then for any $j=1,\cdots,N$, 
 there exists a constant $c_j>0$ and a compact interval $\Delta_j$ about $x$ such that
$$
E_{\lambda_j}(\Delta_j)[\lambda_j,i\mathbf{a}_j]E_{\lambda_j}{(\Delta_j)}
\geq c_jE_{\lambda_j}(\Delta_j).
$$
But, for any Borel set $J$, one has  $E_M(J)=\oplus _{1\leq j\leq N}E_{\lambda_j}(J)$ and 
$[M,i\tilde{A}]=\oplus_{1\leq j\leq N}[\lambda_j,i\mathbf{a}_j]$,
Hence, we immediately  conclude that, for $\Delta=\cap_{j=1}^N\Delta_j$ and $c_j=\min_{1\leq j\leq N}c_j$ we have,
 $$
E_{M}(\Delta)[M,i\tilde{A}]E_{M}{(\Delta)}
\geq cE_{M}(\Delta),
$$
that is, $x\in\mu^{\tilde{A}}(M)$. 
The proof is complete.
%%%%%%%%%%%%%%%%%%%%%%%%%%%%%%%%%%%%%%%%%%%%%%%%%%%%%%%%%%%%%%%%%%%%%%%%%%%%%%%%%%%%%%%%%%%%%%%%%ù
%%%%%%%%%%%%%%%%%%%%%%%%%%%%%%%%%%%%%%%%%%  THEOREM 4.1
%%%%%%%%%%%%%%%%%%%%%%%%%%%%%%%%%%%%%%%%
%%%%%%%%%%%%%%%%%%%%%%%%%%%%%%%%%%%%%%%%%%%%%%%%%%%%%%%%%%%%%%%%%%%%%%%%%%%%%%%%%%%%%%%%%%%%%%%%%

 %%%%%%%%%%%%%%%%%%%%%%%%%%%%%%%%%%%%%%%%%%%%%%%%%%%
According to the Proposition \ref{prop3.1},  $H_0=\mathcal{U}^{-1}M\mathcal{U}$. Define the operator $\mathbb{A}$ by 
\begin{equation}\label{oct}
\mathbb{A}=\mathcal{U}^{-1}\widetilde{A}\mathcal{U}.
\end{equation}
By a direct application of the previous proposition, we get
%%%%%%%%%%%%%%%%%%%%%%%%%%%%%%%%%%%%%%%%%%%%%%%%%%%%%%%%%%%%%%%%%%%%%%%%%%%%%%%%%%%%%%%%%%%%%%%%%ù
%%%%%%%%%%%%%%%%%%%%%%%%%%%%%%%%%%%%%%%%%%  THEOREM 4.1
%%%%%%%%%%%%%%%%%%%%%%%%%%%%%%%%%%%%%%%%
%%%%%%%%%%%%%%%%%%%%%%%%%%%%%%%%%%%%%%%%%%%%%%%%%%%%%%%%%%%%%%%%%%%%%%%%%%%%%%%%%%%%%%%%%%%%%%%%%
\begin{corollary}
\label{moureH_0} The operator $H_{0}$ is of class $C^{\infty }(\mathbb{A})$ and $\mathbb{A}$ is
locally strictly conjugate to $H_{0}$ on $\mathbb{R}\setminus\kappa(H_0)$, i.e.
$
\mu^{\mathbb{A}}(H_0)=\mathbb{R}\setminus\kappa(H_0).
$
\end{corollary}
In particular, the spectrum of $H_0$ is purely absolutely continuous on $\mathbb{R}\setminus\kappa(H_0)$.
In fact,  combining this corollary, Lemma \ref{regularitŽ}, and the results of Section  \ref{rappel}, we immediately get,
%%%%%%%%%%%%%%%%%%%%%%%%%%%%%%%%%%%
%%%%%%%%%%%%%%   LAP   %%%%%%%%%%%%%%%%%
\begin{theorem} \label{rbvcontinuity} 
(i) The  limits 
$
 ({H_0}-x\mp i0)^{-1}:=\lim_{\mu\rightarrow0+}({H_0}-x\mp i\mu)^{-1}
$
exist locally uniformly on  $\mathbb{R}\setminus\kappa(H_0)$ for the weak* topology of 
$B({\mathcal{K},\mathcal{K}^*})$.\\
(ii) For any $s>1/2$, the maps 
$
 x\mapsto\langle {\textbf{N}}\rangle^{-s}({H}_0-x\mp i0)^{-1}\langle {\textbf{N}}\rangle^{-s}\in B({\mathcal   H})
 $
are locally of class $\Lambda^{s-\frac{1}{2}}$.  \\
(iii) For any  $\varphi\in C^\infty_0(\mathbb{R}\setminus\kappa(H_0))$ and  $\sigma,\varepsilon>0$.
$$
||\langle {\textbf{N}}\rangle^{-\sigma}e^{-i{H_0}t}\varphi(H_0)\langle {\textbf{N}}\rangle^{-\sigma}||\leq C(1+|t|)^{-\sigma+\varepsilon}.
$$
\end{theorem}
%%%%%%%%%%%%%%%%%%%%%%%%%%%%%%%%%%%%%%%%%%%%%%%%%%%%%%%%%%%%%%
%%%%%%%%%%%%%%%%%%%%%%%%%%%%%%%%%%%%%%%%%%%%%%%%%%%%%%%%%%%%%%
%%%%%%%%%%%%%%%%%%%%%%%%%%%%%%%%%%%%%%%%%%%%%%%%%%%%%%%%%%%%%%%%%%%%%%%%%
%%%%%%%%%%%%%%%%%%%%%%%%%%%%%%%%%%%%%%%%%%%%%%%%%%%%%%%%%%%%%%%%%%%%%%%%%
\protect\setcounter{equation}{0}
%%%%%%%%%%%%%%%%%%%%%%%%%%%%%%%%%%%%%%%%%%%%%%%%%%%%%%%%%%%%%%%%
%%%%%%%%%%%%%%%%%%%%%%%%%%%%%%%%%%%%%%%%%%%%%%%%%%%%%%%%%%%%%%%%
\section{Mourre estimate for $H$}\label{perturbŽ}
%%%%%%%%%%%%%%%%%%%%%%%%%%%%%%%%%%%%%%%%%%%%%%%%%%%%%%%%%%%%%%%%
%%%%%%%%%%%%%%%%%%%%%%%%%%%%%%%%%%%%%%%%%%%%%%%%%%%%%%%%%%%%%%%%
%
%
In this section we extend our analysis to the operator $H$. 
For we  start by recalling that
$
H=H_0+V
$
where $V$ is the difference operator acting in $\mathcal H$ by
\begin{equation}\label{perturbation}
(V\psi)_n=(A_{n-1}-A)^*\psi_{n-1}+(B_n-B)\psi_n+(A_{n}-A)\psi_{n+1}, ~\mbox{
for all } n\in\mathbb{Z}.
\end{equation}
 According to  (\ref{Jacobi2}),   $V$ is a compact operator in $\mathcal{H}$. 
In particular, 
$$
\sigma_{ess}(H)=\sigma_{ess}(H_0)=\cup_{j=1}^{j=N}\lambda_j([-\pi,\pi])
$$
The main result of this section is the  extension, up to a compact operator, of
 the Mourre estimate obtained in the last section for $H_0$ to 
 $H$. In fact we will prove the following more general result,
\begin{theorem}\label{MEFORH} 
Let $V$ be a self-adjoint compact operator on $\mathcal{H}$ (not necessarily of the form (\ref{perturbation})) such that $V<\textbf{N}>$ is also compact. Then
 the self-adjoint operator $H=H_0+V$  is of class ${C}^{1}(\mathbb{A})$ and $\mathbb{A}$ is locally conjugate to $H$ on $\mathbb{R}\setminus
\kappa (H_0)$, i.e.
$
\tilde{\mu}^{\mathbb{A}}(H)=\mu^{\mathbb{A}}(H_0)=\mathbb{R}\setminus\kappa(H_0).
$
In particular,  outside the finite set $\kappa(H_0)$ the eigenvalues of $H_0+V$ are all finitely degenerate and 
 their possible accumulation points are included in $\kappa(H_0)$. 
\end{theorem}
 \textbf{Proof of Theorem \ref{MEFORH}}
According to the last section, the operator
%%%%%%%%%%%%%%%%%%%%%%%%%%%%%%%%%%%%%%%%%%%%%%%%%%%
$H_0$ is of class $C^\infty(\mathbb{A})$ with
$
\mathbb{A}=\mathcal{U}^{-1}\widetilde{A}\mathcal{U},
$
and $\widetilde{A}$ is defined by (\ref{ocdef1}) and (\ref{defdeatld}). So to prove Theorem \ref{MEFORH} we first show
that $V$ is of class  $C^{1}(\mathbb{A})$ and that $[V,i\mathbb{A}]$ is a compact operator in $\mathcal{H}$.

 According to  Lemma \ref{regularitŽ} below, the operator 
$\langle {\textbf{N}}\rangle^{-1}\mathbb{A}$ is bounded in $\mathcal H$. Hence, 
$
V\mathbb{A}=V\langle {\textbf{N}}\rangle (\langle {\textbf{N}}\rangle^{-1}\mathbb{A})
$
is a compact operator in $\mathcal{H}$. By adjonction, $\mathbb{A}V$ is also a compact operator.
Thus, the commutator $[V,i\mathbb{A}]$ is a compact operator in $\mathcal{H}$.
Moreover, for any Borel set $J$,  $E_H(J)-E_{H_0}(J)$ is a compact operator.
Hence, for any $J\subset\mu^{\mathbb{A}}(H_0)$,
\begin{eqnarray*}
E_H(J)[H,i\mathbb{A}]E_H(J)&=&E_H(J)[H_0,i\mathbb{A}]E_H(J)+E_H(J)[V,i\mathbb{A}]E_H(J)\\
&=&E_{H_0}(J)[H_0,i\mathbb{A}]E_{H_0}(J)+K\\
&\geq&cE_{H_0}(J)+K\\
&\geq&cE_{H}(J)+K'
\end{eqnarray*}
for some $c>0$ and compact operators $K,K'$.  The proof is finished.

\begin{lemma}
\label{regularitŽ} For every integer $m>0$ the operators $\langle {\textbf{N}}\rangle^{-m}\mathbb{A}^{m}$ and
$\mathbb{A}^{m}\langle {\textbf{N}}\rangle^{-m} $ are bounded in $\mathcal{H}$.
\end{lemma}
\noindent \textbf{Proof:}  The proof can be done by induction. We only illustrate the idea for $m=2$ and deal
with $\langle {\textbf{N}}\rangle^{-2} \mathbb{A}^2$ (the remaining  cases are similar). Let us show that this operator is bounded
in $\mathcal{H}$. 

We check that $\mathbb{A}^2 = \mathcal{U}^{-1}\widetilde{A}^2 \mathcal{U}$ is given by 
$$
(\mathbb{A}^2\psi)_n = \sum_{j=1}^{N}\int_{-\pi}^{\pi} \left[G_j^0 (p) \frac{d^2
f_{j}}{dp^2}(p)+ + G_j^1 (p) \frac{d f_{j}}{dp}(p) + G_j^2 (p) f_{j} (p)%
\right] W_{j}(p) e^{-i pn}\, dp, 
$$
where $f_j$ is defined by (\ref{deftrsfF}), and $G_j^0,$ $G_j^1$ and
$G_j^2$
are smooth functions on $[{-\pi},\pi]$.
Integrations by parts show that $(\langle {\textbf{N}}\rangle^{-2} \mathbb{A}^2 \psi)_n$
is a finite sum of terms of the form 
$$
\langle {n}\rangle^{-2} (an^2 +bn +c)\int_{-\pi}^{\pi} f_j (p) W(p) e^{-ipn}\, dp, 
$$
where $a,b$ and $c$ are constants, and $W$ is a smooth  function. Since, on the one hand,
$%
\langle {n}\rangle^{-2} (an^2 +bn +c)$ is bounded, and using, on the other hand,
Plancherel's theorem: 
\begin{eqnarray*}
||\Big(\int_{-\pi}^{\pi} f_j (p) W(p) e^{-ip n}\, dp\Big)_n||_{l^2(\mathbb{Z})}^2&=& ||f_j (p) W(p)||_{L^2([-\pi,\pi])}^2 \\
&\leq& C ||f_j (p)||_{L^2([-\pi,\pi])}^2 \\
&\leq& C ||\psi||_{\mathcal{H}}^2,
\end{eqnarray*}
we obtain $||\langle {\textbf{N}}\rangle^{-2} \mathbb{A}^2 \psi||^2_{\mathcal{H}} \leq C'
||\psi||_{\mathcal{H}}^2$,  for some constants $C,C'>0$.

\textbf{Proof of Theorem \ref{theorem0} }  First, the assumption (\ref{HP1}) implies  that 
$
V\langle {\textbf{N}}\rangle 
$
is a compact operator in $\mathcal{H}$. Hence Theorem \ref{theorem0}  is an immediate consequence of Theorem \ref{MEFORH} and  Theorem \ref{virial}.
\textbf{Proof of Theorem \ref{theorem1} } 
Since condition (\ref{HP}) implies  (\ref{HP1}), $\tilde{\mu}^{\mathbb{A}}(H)=\mathbb{R}\setminus\kappa(H_0)$.
So, according to  Theorem \ref{lap} and Lemma \ref{regularitŽ}, we  only have to prove that $H$ is of class
$\mathcal{C}^{1,1}(\mathbb{A})$. As in the last proof, it is enough to show that
 $V$ is of class  $\mathcal{C}^{1,1}(\mathbb{A})$.  But this is true according to the appendix of \cite{BS1}.
 The proof of Theorem \ref{theorem1} is complete.
 
\textbf{Proof of  Theorem \ref{theorem11} }  As in the precedent proof, according to  Theorem \ref{lap1},  it is enough
to show that the condition (\ref{HP'}) implies  that $V$ is of class  $\mathcal{C}^{1+\theta}(\mathbb{A})$.  
But  (\ref{HP'}) ensures that
the operator $\langle {\textbf{N}}\rangle^{1+\theta }V$ is bounded in $B(\mathcal{H})$. The preceding lemma shows that
 $
 \langle {\mathbb{A}}\rangle^{1+\theta }V=(\langle {\mathbb{A}}\rangle^{1+\theta }\langle {\textbf{N}}\rangle^{-(1+\theta) } )\langle {\textbf{N}}\rangle^{1+\theta }V
 $
is  bounded too.  
The proof of Theorem \ref{theorem11} is finished by the point (ii) of the remark given just after Definition \ref{defreg}.

\textbf{Remark } Here also it should  be clear that Theorem \ref{theorem1} and Theorem \ref{theorem11}  are valid 
for $H=H_0+V$ where the perturbation $V$ is any bounded operator in $\mathcal{H}$ that satisfies one of the following:\\
(i)  $\langle {\textbf{N}}\rangle^{1+\theta }V\in B(\mathcal{H})$ for some $\theta>0$, or\\
(ii) $V$ is a compact operator on $\mathcal{H}$ such that $\langle {\textbf{N}}\rangle^{\theta }[V,\mathbb{A}]\in B(\mathcal{H})$.
%%%%%%%%%%%%%%%%%%%%%%%%%%%%%%%%%%%%%%%%%%
%%%%%%%%%%%%%%%%%%%%%%%%%%%%%%%%%%%%%%%%%%
%%%%%%%%%%%%%%%%%%%%%%%%%%%%%%%%%%%%%%%%%%
%%%%%%%%%%%%    Le cas o $A$ est positive dŽfinie         %%%%%%%%%
%%%%%%%%%%%%%%%%%%%%%%%%%%%%%%%%%%%%%%%%%%
%%%%%%%%%%%%%%%%%%%%%%%%%%%%%%%%%%%%%%%%%%
%%%%%%%%%%%%%%%%%%%%%%%%%%%%%%%%%%%%%%%%%%
\section{The case where $A$ is positive definite }\label{positivedefinie}
%%%%%%%%%%%%%%%%%%%%%%%%%%%%%%%%%%%%%%%%%%
%%%%%%%%%%%%%%%%%%%%%%%%%%%%%%%%%%%%%%%%%%
%%%%%%%%%%%%%%%%%%%%%%%%%%%%%%%%%%%%%%%%%%

In this section  we focus on the case where $H$ is the block Jacobi operator  defined by (\ref{Jacobi1})-(\ref{Jacobi2}) such that $A$ is  positive definite.
Then,  thanks to the Corollary \ref{cas-positive}, we have:
\begin{enumerate}
\item for  all  $ j=1,\cdots,N$, we have
 $$
\lambda_j'(p)=-2<AW_j(p),W_j(p)>_{\mathbb{C}^N}\sin p,\mbox{ for all }p\in[-\pi,\pi].
$$
\item  the eigenvalues $\lambda_1(p),\cdots,\lambda_N(p)$  are even on $[-\pi,\pi]$ and decreasing on $[0,\pi]$.
In particular, $\Sigma_j=[\lambda_j(\pi),\lambda_j(0)]$ and  $\kappa(\lambda_j)=\{\lambda_j(\pi),\lambda_j(0)\}$.
\item the spectrum of $H_0$ is purely absolutely continuous and 
$$
\sigma(H_0)=\cup_{1\leq j\leq N}[\lambda_j(\pi),\lambda_j(0)].
$$
\end{enumerate}
Moreover,  according to the last section, the operator
%%%%%%%%%%%%%%%%%%%%%%%%%%%%%%%%%%%%%%%%%%%%%%%%%%%
$H_0$ is of class $C^\infty(\mathbb{A})$ and $\mu^{\mathbb{A}}(H_0)=\mathbb{R}\setminus\kappa(H_0)$,  with
$
\mathbb{A}=\mathcal{U}^{-1}\widetilde{A}\mathcal{U},
$
and $\widetilde{A}$ is defined by (\ref{ocdef1}) and (\ref{defdeatld}).
The advantage now is that, one may choose 
$$
F_j(p)=-\sin p, j=1,\cdots,N,
$$
in the definition of the operator $\mathbb{A}$.
Indeed,
$$
F_j(p)\lambda_j'(p)=(-\sin p)\lambda_j'(p)=2<AW_j(p),W_j(p)>_{\mathbb{C}^N}(\sin p)^2,
$$
which is strictly positive on $]0,\pi[$. 
But then,  with this choice, one may easily show that 
\begin{equation}\label{opcon3}
\mathbb{A}=\mathbb{D}+L
\end{equation}
for some bounded operator $L$ in  $\mathcal{H}$ and where $\mathbb{D}$ is defined in  $\mathcal{H}$ by,
\begin{equation}\label{opcon2}
(\mathbb{D}\psi)_n=i(n+\frac{1}{2})\psi_{n+1}-i(n-\frac{1}{2})\psi_{n-1}
\end{equation}
This allows us to handle the larger class of perturbations $V$ satisfying
\begin{equation}\label{HP2}
\lim_{|n|\rightarrow\infty}|n|(\|A_{n+1}-A_n\|+\|B_{n+1}-B_n\|)=0.
\end{equation}
Remark that the condition (\ref{HP1}) implies (\ref{HP2}) but not the inverse, think to the example 
$A_n=A+\frac{1}{\ln(1+|n|)}A'$, with any matrix $A'$. So we get the following more general result in comparison with Theorem \ref{theorem0}:
\begin{theorem}\label{thm5.1}
If (\ref{HP2}) holds then $H$ is of class ${C}^{1}(\mathbb{A})$ and 
$\tilde{\mu}^{\mathbb{A}}(H)=\mathbb{R}\setminus\{\alpha_j,\beta_j\}_{1\leq j\leq N}$. In particular, 
 the possible eigenvalues of $H$ in  $\mathbb{R}\setminus\{\alpha_j,\beta_j\}_{1\leq j\leq N}$
are all finitely degenerate and cannot accumulate outside $\{\alpha_j,\beta_j\}_{1\leq j\leq N}$.
\end{theorem}

\textbf{Proof }
The only point we have to verify is the compactness of  $[V,i\mathbb{A}]$ in $\mathcal{H}$
under the conditions  (\ref{Jacobi2}) and  (\ref{HP2}).
Indeed,
$$
[V,i\mathbb{A}]=[V,i\mathbb{D}]+[V,iL].
$$
But $[V,iL]=i(VL-LV)$ is clearly a compact operator since each term is
 ($V$ being compact and $L$ bounded). 
Moreover, a direct computation  shows that the only nonzero matrix coefficients 
of the commutator  $K=[V,i\mathbb{D}]$ are:\\
$
K_{n,n+2}=K_{n+2,n}^*=(n+\frac{1}{2})(A_{n+1}-A_n)-(A_{n+1}-A);$\\
$
 K_{n,n+1}=K_{n+1,n}^*=(n+\frac{1}{2})(B_{n+1}-B_n);$\\
 $ K_{n,n}=(n+\frac{1}{2})[(A_n-A_{n-1})+(A_n-A_{n-1})^*]+(A_{n-1}-A)+(A_{n-1}-A)^*.$\\
 All these coefficients tend to zero at infinity, thanks to  (\ref{Jacobi2}) and  (\ref{HP2}).
 Hence, $[V,i\mathbb{D}]$ is a compact operator and the proof is complete.
  
 \textbf{Remark: } In addition of (\ref{Jacobi2}) suppose that 
 \begin{equation}\label{HP3'}
\int_1^\infty \sup_{r<|n|<2r}\Big(|n|\big(\|A_n-A_{n-1}\|+\|B_n-B_{n-1}\|\big)+\|A_n-A\|\Big) \frac{dr}{r}<\infty.
\end{equation}
Then according to the appendix of \cite{BS1} we get $H$ is of class $\mathcal{C}^{1,1}(\mathbb{A})$ (in fact we have slightly more than that but we will not go in these details here).
 Remark that  (\ref{HP3'}) follows if, for some $\theta>0$, we have
\begin{equation}\label{HP3}
\sup_{n\in\mathbb{Z}}|n|^{\theta}\Big(|n|\big(\|A_n-A_{n-1}\|+\|B_n-B_{n-1}\|\big)+\|A_n-A\|\Big)<\infty.
\end{equation}
In this case it is clear that $\langle \textbf{N}\rangle^\theta[V,i\mathbb{A}]$ is bounded.
Hence  $\langle \mathbb{A}\rangle^\theta[V,i\mathbb{A}]= (\langle \mathbb{A}\rangle^\theta\langle \textbf{N}\rangle^{-\theta})\langle \textbf{N}\rangle^\theta[V,i\mathbb{A}]$
 is  bounded too. Consequently,  $H$ is of class $\mathcal{C}^{1+\theta}(\mathbb{\mathbb{A}})$ so that, we have
 \begin{theorem}\label{thm5.2} Assume (\ref{Jacobi2}). 
Theorems \ref{theorem1}, respectively Theorem \ref{theorem11}, remain valid 
 for this model if we replace the condition (\ref{HP}), respectively (\ref{HP'}),  by  (\ref{HP3'}), respectively (\ref{HP3}).
\end{theorem}
%%%%%%%%%%%%%%%%%%%%%%%%%%%%%%%%%%%%%
%%%%%%%%%%%%%%%%%%%%%%%%%%%%%%%%%%%%%
%%%%%%%%%%%%    OpŽrateurs de Schrodinger    %%%%%%%%%
%%%%%%%%%%%%     avec potentiel anisotrope      %%%%%%%%%
%
%
\begin{example} Let us consider  the Jacobi operator $J$ acting in  $l^2(\mathbb{Z},\mathbb{C})$   by
$$
(J\psi)_n=a_{n-1}\psi_{n-1}+b_n\psi_n+a_{n}\psi_{n+1},
~~\mbox{ for all }n\in\mathbb{Z}.
$$
Here-above $a_n>0$ and $b_n\in\mathbb{R}$.  Assume that there exist $a^\pm>0,b^\pm\in\mathbb{R}$ such that,
\begin{equation}\label{aabb}
\lim_{n\rightarrow\pm\infty}a_n=a^\pm~~\mbox{and }~~\lim_{n\rightarrow\pm\infty}b_n=b^\pm.
\end{equation}
This model has been studied in \cite{BS2} where the authors used some ideas of three body problem to establish their Mourre estimate. 
Here our present approach applies in straightforward
way. Indeed, Let $U:l^2(\mathbb{Z})\rightarrow l^2(\mathbb{N},\mathbb{C}^2)$
defined by
$
(U\gy)_n=(\gy_n,\gy_{-n+1}).
$
It is clear, see \cite{B}, that $U$ is unitary operator and 
$UJU^{-1}=H$ 
where  $H$ is the block Jacobi operator realized by
$$
 A_n= \left(%
\begin{array}{ccc}
a_n&0\\
0&a_{-n}\\
\end{array}\right),~\mbox{$n\geq1$};
B_n=\left(%
\begin{array}{ccc}
b_n&0\\
0&b_{-n+1}\\
\end{array}\right),~\mbox{$n\geq2$}~\mbox{ and }~
B_1:=\left(%
\begin{array}{ccc}
b_1&a_0\\
a_0&b_{0}\\
\end{array}\right).
$$
In particular,  one may deduce directly by applying the last discussion that $\sigma_{ess}(J)=[b^--2a^-,b^-+2a^-]\cup[b^+-2a^+,b^++2a^+]$. Moreover, if 
$$
\lim_{n\rightarrow\pm\infty}|n|^{\theta}\Big(|n|\big(|a_n-a_{n+1}|+|b_n-b_{n+1}|\big)+|a_n-a^\pm|)=0,\quad\mbox{for some $\theta\geq0$}
$$
holds then  $\sigma_p(J)$  has no accumulation points outside $\kappa(J)=\{b^-\pm2a^-,b^+\pm2a^+\}$. Moreover, if
 $\theta>0$ then  $\sigma_{sc}(J)=\emptyset$.
 \end{example}
 % %%%%%%%%%%%%%%%%%%%%%%%%%%%%%%%%%%%%%%%%%%
%%%%%%%%%%%%%%%%%%%%%%%%%%%%%%%%%%%%%%%%%%
%%%%%%%%%%%%%%%%%%%%%%%%%%%%%%%%%%%%%%%%%%
\section{Difference operators on cylindrical domains }\label{waves}
%%%%%%%%%%%%%%%%%%%%%%%%%%%%%%%%%%%%%%%%%%
%%%%%%%%%%%%%%%%%%%%%%%%%%%%%%%%%%%%%%%%%%
%
 Consider the configuration space $X=\mathbb{Z}\times \{1,\cdots,N\}$,  for some
  integer $N$, and  the difference operator $J$ acting in $l^2(X)$   by
\begin{equation}\label{Jacobicylinder}
(J\psi)_{n,k}=a_{n-1,k}\psi_{n-1,k}+b_{n,k}\psi_{n,k}+a_{n,k}\psi_{n+1,k}+
c_{n,k-1}\psi_{n,k-1}+c_{n,k}\psi_{n,k+1}
\end{equation}
initialized by $\psi_{n,0}=\psi_{n,N+1}=0$ and where $a_{n,k},c_{n,k}$ are positive while $b_{n,k}$ are real numbers.
\begin{theorem}\label{dwp}
Assume that, there is some positive numbers $a_1,\cdots,a_N, c_1,\cdots,c_{N-1}$ and real numbers $b_1,\cdots,b_N$ such that for all $j=1,\cdots,N $,
$$
\lim_{n\rightarrow\pm\infty}|a_{n,j}-a_j|+|b_{n,j}-b_j|+|c_{n,j}-c_j|=0~,\quad (\mbox{with the convention $c_{n,N}=c_{N}=0$}).
$$
Then $J$ is a bounded symmetric operator whose essential spectrum is of the form 
$\cup_{j=1}^N[\alpha_j,\beta_j]$ for some numbers
$\alpha_j<\beta_j$. If in addition, we have
$$
\lim_{n\rightarrow\pm\infty}|n|(|a_{n+1,j}-a_{n,j}|+|b_{n+1,j}-b_{n,j}|+|c_{n+1,j}-c_{n,j}|)=0, 
$$
then there is  a self-adjoint operator $\mathbb{A}'$ such that $\tilde{\mu}^{\mathbb{A}'}(J)=\mathbb{R}\setminus\kappa(J)$ with 
$\kappa(J)=\cup_{1\leq j\leq N}\{\alpha_j,\beta_j\}$. In particular,  the possible eigenvalues of $J$ are all finitely degenerate and cannot accumulate  outside 
$\kappa(J)$.
Finally,  $\sigma_{sc}(J)=\emptyset$  and conclusions of Theorems \ref{theorem1} and
\ref{theorem11} hold for $J$ if, for some $\theta>0$, we have
$$
\sup_{n}|n|^{\theta}{\Big(|n|\big(|a_{n+1,j}-a_{n,j}|+|b_{n+1,j}-b_{n,j}|+|c_{n+1,j}-c_{n,j}|\big)+|a_{n,j}-a_j|\Big)}<\infty.
$$ 
\end{theorem}
Before to show this result let us fix some notations that we will use in the sequel of this paper. 
Let $E_{ij}$  be the elementary matrix whose only non zero entry is 1 placed at the intersection of
the $i^{th}$ ligne and the $j^{th}$ column  and
\begin{eqnarray}\label{diagonal-tridiagonal}
 \mbox{diag}(d_1,\cdots,d_N)&=& \left(%
\begin{array}{cccc}
d_{1}&0&\dots&0\\
0&d_{2}&0&\dots\\
\ddots&\ddots&\ddots&\ddots\\
0&\dots&0&d_{N}\\
\end{array}\right)
\end{eqnarray}
\begin{eqnarray}
\mbox{tridiag}(\{\alpha_i\}_{1\leq i\leq N-1},\{\beta_i\}_{1\leq i\leq N})&=&\left(%
\begin{array}{ccccccc}
\beta_{1}&\alpha_{1}&0&\dots&0\\
\alpha_{1}&\beta_{2}&\alpha_{2}&0&\dots\\
0&\alpha_{2}&\beta_{3}&\alpha_{2}&\dots\\
\ddots&\ddots&\ddots&\ddots&\ddots\\
0&\dots&0&\alpha_{N-1}&\beta_{N}\\
\end{array}\right)~~
\end{eqnarray}

\textbf{Proof of Theorem \ref{dwp}  }
 Let  $U:l^2(\mathbb{Z}\times \{1,\cdots,N\})
\rightarrow l^2(\mathbb{Z},\mathbb{C}^N)$
defined by
$
(U\gy)_n= (
\gy_{n,1},\cdots,
\gy_{n,N}).
$
It is clear that $U$ is unitary operator and  
$
UJU^{-1}=H(\{A_n\},\{B_n\})
$ where 
\begin{eqnarray}\label{rho}
 A_n=\mbox{diag}(a_{n,1},\cdots,a_{n,N})
 ~\mbox{ , }~B_n=\mbox{tridiag}(\{c_{n,i}\}_{1\leq i\leq N-1},\{b_{n,i}\}_{1\leq i\leq N})
\end{eqnarray}
forall $n\in\mathbb{Z}$. By hypothesis  the matrices
 $A_n$ and $B_n$  satisfy (\ref{Jacobi2})  with
\begin{eqnarray}\label{8.12}
A= \mbox{diag}(a_{1},\cdots,a_{N})\quad\mbox{and}\quad 
B=\mbox{tridiag}(\{c_{i}\}_{1\leq i\leq N-1},\{b_{i}\}_{1\leq i\leq N})
\end{eqnarray}
Since $A$ is positive definite,  the results of the last section apply directly to the operator $H(\{A_n\},\{B_n\})$. 
In addition, the conjugate operator of $J$ is given by $\mathbb{A}'=U^{-1}\mathbb{A}U$ where $\mathbb{A}$ is defined by 
(\ref{opcon3}) and (\ref{opcon2}).
The proof is finished.
 \begin{example}
 To have  explicit formulae, assume in (\ref{8.12}) that
 $$
 a_i=a>0, ~~b_i=0~~ \mbox{ and~ $c_i=c>0$.}
 $$
 In this case,  $h(p)$ is $N$ by $N$  Toeplitz matrix whose spectrum is given by
 $$
 \lambda_j(p)=2a\cos p+2c\epsilon_j~~Ê\mbox{ with }~~\epsilon_j=\cos\frac{(N-j+1)\pi}{N+1}~,~ j=1,\cdots,N.
 $$
 Here again, we observe explicitly
 $
 \sigma(H_0)=\cup_{j=1}^N[\lambda_j(\pi),\lambda_j(0)]
 $.
 Remark that $j^{th}$ gap is non trivial, if and only if, 
 $$
L_j=\lambda_{j+1}(\pi)-\lambda_j(0)=2c\Delta_j-4a>0\quad\mbox{where }
\Delta_j=\epsilon_{j+1}-\epsilon_{j}>0.
$$
Then $j^{th}$ gap is non trivial, if and only if,
$
c>{2a}/{\Delta_j}.
$
  \end{example}
\begin{example}\textbf{Discrete wave propagator. }
Consider on $X=\mathbb{Z}\times \{1,\cdots,N\}$ 
the discrete Laplacian $\Delta$ given by
$$
(\Delta\gy)_{n,m}=\gy_{n+1,m}+\gy_{n-1,m}+\gy_{n,m+1}+\gy_{n,m-1}~,~\mbox{for
  all $x=(n,m)\in X$}
$$
with the initial boundary condition
$\gy_{0,m}=\gy_{n,0}=\gy_{n,N+1}=0$.
Let $\rho:X\rightarrow]0,\infty[$ be the local propagation speed. The 
 discrete wave equation is given by
$$
\frac{\partial^2u}{\partial{t}^2}=\rho^2\Delta u.
$$
Clearly, the change of variable $u=\rho v$ transforms the last
equation  to
$$
\frac{\partial^2v}{\partial{t}^2}=\rho\Delta \rho v,
$$
Hence the spectral analysis of the operator  $H=\rho\Delta \rho$ on $\mathcal{H}=l^2(X)$ is of interest for discrete wave
equation. Assume that the cylinder $X$ is stratified, i.e. for some constant $\rho_1,\cdots,\rho_N>0$, 
one has $\rho(n,m)-\rho_m\rightarrow0$ as
$n\rightarrow\infty$. Hence  $H=H(\{A_n\},\{B_n\})$ given in (\ref{rho}) with 
$$
a_{n,m}=\rho(n,m)\rho(n+1,m)~~,~~b_{n,m}=0~\mbox{ and } ~ 
c_{n,m}=\rho(n,m)\rho(n,m+1).
$$
Of course we assume that $\rho_0=\rho_{N+1}=0$.
This case can be clearly studied by Theorem \ref{dwp} since
$$
 a_{n,m}\rightarrow a_m=\rho_m^2~\mbox{ and } 
c_{n,m}\rightarrow c_m=\rho_m\rho_{m+1}~\mbox{ as } |n|\rightarrow\infty.
$$
But we will not state the corresponding result separately. 

For example, if the stratification is simple in the sense that $\rho$ is constant, say  $\rho_m=1$ for all $m$.
 Then the essential spectrum of $H$ has no gaps since the spectral bands $\Sigma_j$ and $\Sigma_{j+1}$ overlap. Indeed, with the notations of example 8.1 we have:
 $$
L_j=\lambda_{j+1}(\pi)-\lambda_j(0)=2\Delta_j-4\leq0\quad\mbox{here }
\Delta_j=\epsilon_{j+1}-\epsilon_{j}>0
$$
\end{example}
%%%%%%%%%%%%%%%%%%%%%%%%%%%%%%%%%%%%%%%%%%%%%%%%%
%%%%%%%%%%%%%%%%%%%%%%%%%%%%%%%%%%%%%%%%%%%%%%%%%
%%%%%%%%%%%%%%%%%%%%%%%%%%%%%%%%%%%%%%%%%%%%%%%%%
%\subsection{ A first model where $A$ is singular }
%%%%%%%%%%%%%%%%%%%%%%%%%%%%%%%%%%%%%%%%%%%%%%%%%
%%%%%%%%%%%%%%%%%%%%%%%%%%%%%%%%%%%%%%%%%%%%%%%%%
%%%%%%%%%%%%%%%%%%%%%%%%%%%%%%%%%%%%%%%%%%%%%%%%%
\protect\setcounter{equation}{0}
\section{A first case where $A$ is singular}\label{app-singular1}
 The simplest situation one may consider is a
naive modification of the example considered in the previous section, namely where $B$ is a tridiagonal matrix and  
$A$ has only one nonzero coefficient on the diagonal.
More specifically, let $a_1,a_2,\cdots,a_N>0$ and $b_1,\cdots,b_N\in\mathbb{R}$ be given and 
$H$ be the block Jacobi operator  defined by (\ref{Jacobi1})  
where $A_n$ and $B_n$  satisfy (\ref{Jacobi2})  with
\begin{eqnarray}
A=a_NE_{N,N}\quad \mbox{and}\quad B=\mbox{tridiag}(\{a_{i}\}_{1\leq i\leq N-1},\{b_{i}\}_{1\leq i\leq N})
\end{eqnarray}
\begin{theorem} \label{singulier1}
\begin{enumerate}
\item There exist $\alpha_1<\beta_1<\alpha_2<\beta_2\cdots<\alpha_N<\beta_N$  such that  the spectrum of $H_0$ is purely absolutely continuous and
$
\sigma(H_0)=\cup_{j=1}^N[\alpha_j,\beta_j]
$
and $\kappa(H_0)=\{\alpha_1,\beta_1,\alpha_2,\beta_2\cdots,\alpha_N,\beta_N\}$. 
\item
If  (\ref{HP1}) holds then  we have a Mourre estimate for $H_0$ 
and $H$  on $\mathbb{R}\setminus\kappa(H_0)$, i.e.
$
\mu^{\mathbb{A}}(H_0)=\tilde{\mu}^{\mathbb{A}}(H)=\mathbb{R}\setminus\kappa(H_0).
$ In particular,  the eigenvalues of $H$ are all finitely degenerate and cannot accumulate outside $\kappa(H_0)$.
\item If  (\ref{HP})  holds then the conclusions of Theorems \ref{theorem1} and \ref{theorem11} are valid.
\end{enumerate}
\end{theorem}
Mention here that the spectrum of $H_0$ has a band/gap structure, that is between each two successive bands $\Sigma_j=[\alpha_j,\beta_j]$ and
$\Sigma_{j+1}=[\alpha_{j+1},\beta_{j+1}]$ there is a non trivial gap of length $L_j=\alpha_{j+1}-\beta_{j}>0$.
To show this result we will study  the eigenvalues of $h(p)$  given here   by
$$
h(p)=\mbox{tridiag}(\{a_{i}\}_{1\leq i\leq N-1},\{b_1,\cdots,b_{N-1},b_N+2a_N\cos p\})
$$
Since $h(-p)=h(p)$, it  is enough to restrict ourselves to $p\in[0,\pi]$. 
Theorem \ref{singulier1} is an  immediate consequence of   the following:
\begin{theorem}\label{singulier0}
Let $B_{N-1}=\mbox{tridiag}(\{a_{i}\}_{1\leq i\leq N-2},\{b_{i}\}_{1\leq i\leq N-1})$ and $\mu_1<\mu_1\cdots<\mu_{N-1}$ its eigenvalues.
The eigenvalues $\lambda_1(p),\cdots,\lambda_N(p)$  are simple,  decreasing on $[0,\pi]$ and can be chosen so that, for all $p\in[0,\pi]$,
\begin{equation}\label{naive}
\lambda_1(p)<\mu_1<\lambda_2(p)\cdots<\mu_{N-1}<\lambda_N(p).
\end{equation} 
In particular, $\lambda_j([0,\pi])=[\alpha_j,\beta_j]$ with $\alpha_j=\lambda_j(\pi)$ and  $\beta_j=\lambda_j(0)$.  
Moreover, $\kappa(\lambda_j)=\{\alpha_j,\beta_i\}$.
\end{theorem}
\textbf{Proof }
Let $p\in[0,\pi]$.  The matrix $h(p)$ is  symmetric tridiagonal with positive off-diagonal elements and
 its characteristic polynomial  satisfies
$$
D(x,p)=\det(h(p)-x)=(b_N+2a_N\cos p-x)D_{N-1}(x)-a_{N-1}^2D_{N-2}(x),
$$
where $D_{N-j}(x)=\det(B_{N-j}-x)$ and $B_{N-j}$ is obtained from $h(p)$, and so from $B$, by eliminating the last $j$ rows and columns.
Moreover,  it is known  that $D(x,p)$ has $N$ distincts real roots that are separated by the $N-1$ distinct roots of $D_{N-1}(x)$.
So (\ref{naive}) is proved. Moreover,  direct calculation shows that
$$
\frac{\partial D}{\partial p}(\lambda_j(p),p)=-2a_N(\sin{p}) D_{N-1}(\lambda_j(p))\not=0~,~\mbox{ for all $p\in]0,\pi[$}.
$$
 Thus
by differentiating the identity $D(\lambda_j(p),p)=0$ we get
$$
\lambda'_j(p)\frac{\partial D}{\partial x}(\lambda_j(p),p)+\frac{\partial D}{\partial p}(\lambda_j(p),p)=0
$$
and so $\lambda'_j(p)$  is non zero on $]0,\pi[$. Finally, 
because the matrix $A$ is non negative we have, for 
 all  $j=1,\cdots,N$,
 $$
\lambda_j'(p)=-2<AW_j(p),W_j(p)>_{\mathbb{C}^N}\sin p\leq0, \quad\mbox{for all }p\in[0,\pi].
$$
 Hence $\lambda'_j(p)<0$ for all $p\in]0,\pi[$.
We see that $\kappa(\lambda_i)=\{\lambda_j(0),\lambda_j(\pi)\}$. The proof is finished.

\textbf{Remark: }
Mention that if we put the nonzero coefficient at the first place  of the diagonal instead of the last one then 
exactly the same phenomenon happens. 
In contrast, if we put it in a different place of the diagonal then
some drastic changes may arise. For example, assume that $B$ is the same and  $A=a_NE_{N-1,N-1}$. Hence
$$
 h(p)=\mbox{tridiag}(\{a_{i}\}_{1\leq i\leq N-1},\{b_1,\cdots,b_{N-2},b_{N-1}+2a_N\cos p,b_N\})
$$
In  this case, on one hand
$
D(x,p)=\det(B-x)+2a_{N}\cos{p}(b_N-x)D_{N-2}(x).
$
On the other hand
$
\det(B-x)=(b_N-x)D_{N-1}(x)-a_{N-1}^2D_{N-2}(x).
$
Hence,
$$
D(x,p)=(b_N-x)\Big(D_{N-1}(x)+2a_{N}\cos{p}(b_N-x)D_{N-2}(x)\Big)-a_{N-1}^2D_{N-2}(x).
$$
Now let us choose first $a_1,\cdots a_{N-3}$ and $b_1,\cdots,b_{N-2}$ then  pick $b_N$ 
among the eigenvalues of $B_{N-2}$. So $b_N-x$ divide $D_{N-2}(x)$. Hence, $b_N$ is a root of $D(x,p)$ for all $p$ and therefore
$H_0$ has a degenerate spectral band and $b_N$ is an infinitely degenerate eigenvalue of $H_0$.

%%%%
\begin{example}\label{inverse1} Let  $N=2$ so that $A=a_{2}E_{2,2} ~,~B=\mbox{tridiag}(\{a_1\},\{b_1,b_2\})$
where $a_1,a_2>0 $ and $b_1,b_2\in\mathbb{R}$. It is easy to check that
$\sigma(H_0)=[\lambda_1(\pi),\lambda_1(0)]\cup[\lambda_2(\pi),\lambda_2(0)]$ and that $\lambda_1(0)<b_1<\lambda_2(\pi)$.

Let us consider now the inverse problem of this model. 
More specifically, let $\alpha_1<\beta_1,\alpha_2<\beta_2$ be given and find $a_1,a_2>0 $ and $b_1,b_2\in\mathbb{R}$ such that the corresponding $H_0$ satisfies
$
\sigma(H_0)=[\alpha_1,\beta_1]\cup[\alpha_2,\beta_2].
$
First, if there is a solution then, according to Theorem \ref{singulier1} and Theorem \ref{singulier0},  
$\alpha_1,\alpha_2$ are the eigenvalues of $h(\pi)$  and $\beta_1,\beta_2$ are the eigenvalues of  $h(0)$. Thus
$$
\begin{array}{ccc}
  \alpha_1+\alpha_2=b_1+b_2-2a_2& ,  &\alpha_1\alpha_2=b_1(b_2-2a_2)-a_1^2   \\
 \beta_1+\beta_2=b_1+b_2+2a_2 &,  & \beta_1\beta_2=b_1(b_2+2a_2)-a_1^2   
\end{array}
$$
Hence $4a_2=(\beta_1-\alpha_1)+(\beta_2-\alpha_2)$ and $4a_2b_1=\beta_1\beta_2-\alpha_1\alpha_2$ gives $b_1$ while
$2(b_2+b_1)=\alpha_1+\alpha_2+\beta_1+\beta_2$ leads to $b_2$. Finally we get $2a_1^2=2b_1b_2-\alpha_1\alpha_2-\beta_1\beta_2$ which gives $a_1$.
Of course we have to make sur that $2b_1b_2-\alpha_1\alpha_2-\beta_1\beta_2>0$ which holds  if and only if  $\beta_1<\alpha_2$. 
\end{example}
%%%%
%%%%
%%
%
%
%%%%%%%%%%%%%%%%%%%%%%%%%%%%%%%%%%%%%%%%%%%%%%%%%
%%%%%%%%%%%%%%%%%%%%%%%%%%%%%%%%%%%%%%%%%%%%%%%%%
%%%%%%%%%%%%%%%%%%%%%%%%%%%%%%%%%%%%%%%%%%%%%%%%%
%%%%%%%%%%%%%%%%%%%%%%%%%%%%%%%%%%%%%%%%%%%%%%%%%
\protect\setcounter{equation}{0}
\section{A second case where $A$ is singular:  Periodic Jacobi operators}\label{periodique}
In this part we  look at  the case where the only nonzero entry of $A$ is put at the lower left corner and $B$ is a tridiagonal matrix.
It turns out that scalar Jacobi operators with periodic coefficients are covered by this model. More precisely, let $H_0=H_0(A,B)$  defined by (\ref{Jacobi0}) with 
\begin{eqnarray}\label{A,Bperiodique}
A=a_{N}E_{N,1} ~,~B=\mbox{tridiag}(\{a_{i}\}_{1\leq i\leq N-1},\{b_i\}_{1\leq i\leq N}).
\end{eqnarray}
Then we have
\begin{theorem} \label{thm7.1}
\begin{enumerate}
\item There exist $\alpha_1<\beta_1\leq\alpha_2<\beta_2\cdots\leq\alpha_N<\beta_N$  such that  the spectrum of $H_0$ is purely absolutely continuous and
$
\sigma(H_0)=\cup_{j=1}^N[\alpha_j,\beta_j].
$
\item Define the critical set $\kappa(H_0)=\{\alpha_1,\beta_1,\alpha_2,\beta_2\cdots,\alpha_N,\beta_N\}$. 
There exists a conjugate operator $\mathbb{A}$ for $H_0$  i.e.
$
\mu^{\mathbb{A}}(H_0)=\mathbb{R}\setminus\kappa(H_0).
$ 
In particular,  conclusions of Theorem \ref{rbvcontinuity} hold for $H_0$.
\item Let $V$ be a symmetric compact operator. Then the sum $H=H_0+V$  defines a bounded self-adjoint operator in $\mathcal{H}$
 whose essential spectrum coincides with $\sigma(H_0)$. Moreover, if 
 $\langle\textbf{N}\rangle V$ is  compact  then 
 $
\tilde\mu^{\mathbb{A}}(H)=\mathbb{R}\setminus\kappa(H_0).
$ 
  In particular, the eigenvalues of $H$ are all finitely degenerate and cannot accumulate outside $\kappa(H_0)$.
\item If $\langle\textbf{N}\rangle^{1+\theta} V$ is bounded for some $\theta>0$ then the conclusions of Theorems \ref{theorem1} and \ref{theorem11} are valid for $H$.
\end{enumerate}
\end{theorem}
A consequence of this theorem is the following result on Jacobi operators with periodic coefficients,
\begin{corollary}
Let  $J_0$ be the periodic Jacobi operator given in $l^2(\mathbb{Z})$ by
$$
(J_0 \psi)_n=a_n\psi_{n+1}+b_n\psi_n+a_{n-1}\psi_{n-1}
$$
such that $a_n>0$ and  $b_n\in\mathbb{R}$ with $a_{j+N}=a_j$ and $b_{j+N}=b_j$.
Then assertions of Theorem \ref{thm7.1} hold true when we replace $H_0$ by $J_0$. 
\end{corollary}
\textbf{Proof }
 Consider the unitary operator $U:l^2(\mathbb{Z})\rightarrow l^2(\mathbb{Z},\mathbb{C}^N)$ defined by
$(U\psi)_n\!=\!(\psi_{(n-1)N+1},    \cdots,\psi_{nN} )$. 
One may show that 
$
UJ_0U^{-1}=H_0(A,B)
$
where $A,B$ are given by (\ref{A,Bperiodique}). The proof is complete.

Of course spectral properties  of $J_0$ are well known, see \cite{T} for example.
 The point here is to obtain the Mourre estimate for these operators and also to show how our general approach works in this context. 

Here again spectral analysis of   $H_0$  is based on the study  of the eigenvalues of the reduced operators $h(p)$ given here by
\begin{equation}\label{h(p)periodique}
h(p)=\left(
\begin{array}{ccccc}
  b_1&a_1&0&\cdots&a_Ne^{ip}  \\
  a_1&b_2&a_2&0&\cdots\\
  0&\ddots&\ddots&\ddots&\vdots\\
  0&\cdots&a_{N-2}&b_{N-1}&a_{N-1}\\
  a_Ne^{-ip}&\cdots&0&a_{N-1}&  b_N  
\end{array}
\right) ~,\quad \forall p\in[-\pi,\pi] .
\end{equation}
Since $h(-p)$ is the transpose of $h(p)$, it  is enough to restrict ourselves to $p\in[0,\pi]$. 
Theorem \ref{thm7.1} is an  immediate consequence of   the following:
\begin{theorem}\label{thm7.0}
The eigenvalues $\lambda_1(p), \cdots,\lambda_N(p),$ of $h(p)$ are  monotonic on $[0,\pi]$ and can be chosen so that
$\lambda_1(p)<\lambda_2(p)\cdots<\lambda_N(p)$, for all $p\in]0,\pi[$. In fact, 
if $\lambda_j$ is increasing (respectively decreasing) then $\lambda_{j+1}$ is decreasing (respectively increasing).
 Moreover, for $p=0$ or $\pi$ the eigenvalues $\lambda_j(p)$  may be double and
 \begin{eqnarray}\label{eigenvalues}
\lambda_1(\pi)\!<\!\lambda_1(0)\!\leq\!\lambda_2(0)\!<\!\lambda_2(\pi)\!\leq\!\cdots\!\leq\!\lambda_N(\pi)\!<\!\lambda_N(0)
 \mbox{ if $N$ is odd };\quad\quad\\
\lambda_1(0)\!<\!\lambda_1(\pi)\!\leq\!\lambda_2(\pi)\!<\!\lambda_2(0)\!\leq\!\cdots\!\leq\!\lambda_N(\pi)<\!\lambda_N(0)
 \mbox{ if $N$ is even }.\label{eigenvalues1}\quad\quad
\end{eqnarray}
In particular,  there exist $\alpha_1<\beta_1\leq\alpha_2<\beta_2\cdots\alpha_{N-1}<\beta_{N-1}\leq\alpha_N<\beta_N$  such that $\lambda_j([0,\pi])=[\alpha_j,\beta_j]$ and
$\kappa(\lambda_j)=\{\alpha_j,\beta_j\}$.
If for some $1\leq  j\leq N$,  $\beta_j=\alpha_{j+1}$ then such a value
is a double eigenvalue of $h(p)$ for $p=0$ or $\pi$.
\end{theorem}
Notice that  $h(p)$ for $p=0$ and $\pi$ are exactly the periodic and anti-periodic Jacobi matrices
studied in \cite{KvM,vM} where the assertion (\ref{eigenvalues})-(\ref{eigenvalues1}) are proved.
More generally, eigenvectors of $h(p), p\in[0,\pi]$ are nothing but the so-called Bloch solutions studied in the Floquet theory 
developed for periodic Jacobi operators.  In this sense, Theorem \ref{thm7.0}  is well  known with different proofs, see \cite{DKS,KvM,vM,T}. 
We will give however a  proof, mainly for completeness, but also because it is elementary and based on standard linear algebra.
It will be given in few lemmas.
\begin{lemma}\label{lemma1}
For every $p\in]0,\pi[$, the eigenvalues $\lambda_1(p),\cdots,\lambda_N(p)$ of $h(p)$ are all simple,
 while for $p=0,\pi$ they are at most  of multiplicity two. Moreover, they are all monotonic on $(0,\pi)$.
\end{lemma}
\textbf{Proof  } Let  $u=(x_1,\cdots,x_N)$ be an eigenvector of $h(p)$ corresponding to an eigenvalues $\lambda$.
This is equivalent to,
\begin{equation}\label{eigenvaluepbm}
\left\{
\begin{array}{cccccccc}
(b_1-\lambda)x_1  &+&a_1x_2&+ &a_Ne^{ip}x_N&=0 &  \\
a_{i-1}x_{i-1}&+&(b_i-\lambda)x_i&+&  a_ix_{i+1}&=0&, ~~~2\leq i\leq N-1   \\
 a_Ne^{-ip}x _{1}&+&a_{N-1}x_{N-1} &+& (b_N-\lambda)x_N&=0&  
\end{array}
\right.
\end{equation}
We see that for all $2\leq i\leq N-1$, $x_i$ is a linear combination of $x_1$ and $x_N$ so that $\lambda$ is at most of multiplicity two.
 Moreover, if $x_1=0$ then $p=0$ or $p=\pi$. Indeed, if $x_1=0$ then $x_2=-(a_N/a_1)e^{ip}x_N$. Moreover, using simultanousely 
  the equations $j=2,3,\cdots,N-1$, we get that $x_N=c e^{ip}x_N$ for some $c\in\mathbb{R}$.  This means that $x_N=0$ or $e^{ip}\in\mathbb{R}$.
Since $u$ is an eigenvector of $h(p)$,  $u\not=0$ so that $x_N\not=0$ and therefore $e^{ip}\in\mathbb{R}$,or equivalently $p=0$ or $p=\pi$. Similarly one may prove that if $x_N=0$ then $p=0$ or $p=\pi$.

Now let $p\in(0,\pi)$ and assume that $e^{ip}\overline{x_1}x_N$ is real. Then multiplying our system by $\overline{x_1}$ we get from the first equation that $\overline{x_1}x_2$ is real too. Similarly, using simultanousely 
 the equations $j=2,3,\cdots,N-1$, we get $\overline{x_1}x_N$ is real which is impossible. Thus $\Im(e^{ip}\overline{x_1}x_N)\not=0$ and so $\lambda'(p)\not=0$ according to the proposition \ref{deŽrivŽee}.
This completes the proof.
\begin{lemma}\label{lemma2} Let $a=(-1)^{N-1}a_1a_2\cdots a_{N}$. The characteristic polynomial of $h(p)$ is given by,
\begin{eqnarray}\label{D(x,p)}
D(x,p)&:=&\det(h(p)-x)=\Delta(x)+2a\cos p,\\
\Delta(x)&=&\det(B-x)-a_{N}^2\det(\dot{B}-x)
\end{eqnarray}
where 
 $\dot{B}$ is the matrix obtained  from $B$ by removing the first and last rows and columns with the convention $\det(\dot{B}-x)=1$ in the case where $N=2$.
\end{lemma}
\textbf{Proof  } Direct calculation, see however \cite{H} where the formula is already used.

\textbf{Remark } We  emphasize  that, according to (\ref{D(x,p)}),   if one of the $a_j$'s is zero then $D(x,p)=\Delta(x)$ which  is independent of $p$.  In this case,  the spectrum of $h(p)$ is independent of $p$  so that the spectrum of $H_0$ consists of $N$ infinitely degenerate eigenvalues. This provides one with a general class of block Jacobi matrices with this kind of spectra.
\begin{lemma}\label{lemma3}
 The derivative $\Delta'$ of $\Delta$ has $N-1$ roots
 $\mu_1<\mu_2<\cdots<\mu_{N-1}$ and 
the $\lambda_j$'s may be chosen so that, for all $j=1,\cdots,N$, \\
(i) $(-1)^j\Delta'(\lambda_j(p))>0$ for all  $ p\in]0,\pi[$;\\
(ii) $\lambda_j(]0,\pi[)\subset]\mu_{j-1},\mu_j[$ with $\mu_0=-\infty$ and $\mu_N=+\infty$.\\
(iii) For any $j=1,\cdots,N$,  $(-1)^{N+j-1}\lambda_j(p)$ is an increasing function of $p\in[0,\pi]$ and $\kappa(\lambda_j)\subset\{\lambda_j(0),\lambda_j(\pi)\}$. 
\end{lemma}
\textbf{Proof }  (i) Let us denote by $\lambda_1(\pi/2)<\lambda_2(\pi/2)<\dots<\lambda_N(\pi/2)$ the distinct eigenvalues of $h(\pi/2)$ (thanks to the Lemma \ref{lemma1}). 
 They are the roots of $\Delta(x)=D(x,\pi/2)$.  So for any $j=1,\cdots,N-1$ there is
 $\mu_j\in]\lambda_j(\pi/2),\lambda_{j+1}(\pi/2)[$ such that
$
\Delta'(\mu_j)=0,
$
Moreover, since $\Delta$ is a polynomial of the form $\Delta(x)=(-1)^Nx^N+\mbox{lower terms}$,  it is clear that $(-1)^j\Delta'(\lambda_j(\pi/2))>0$ for any $j=1,\cdots N$ which gives us the sign of $\Delta'$ on each $]\mu_{j-1},\mu_j[$.

(ii) Now assume that  (ii) is not true for some $j$, say $j=1$  (the other cases are similar). Then,  by continuity of  $\lambda_1$ and the fact that $\lambda_1(\pi/2)<\mu_1$, $\mu_1$ would belong to $ \lambda_1(]0,\pi[)$.
Therefore $\mu_1=\lambda_1(p_0)$ for some $p_0\in]0,\pi[$. But then for such $p_0$ one has:
$$
D(\mu_1,p_0)=D(\lambda_1(p_0),p_0)=\frac{\partial D}{\partial x}(\lambda_1(p_0),p_0)=0.
$$
Here-above we used the fact that $\frac{\partial D}{\partial x}(x,p)=\Delta'(x)$.
This means that $\mu_1=\lambda_1(p_0)$ is a multiple eigenvalue of $h(p_0)$ which, according to Lemma \ref{lemma1},  implies that $p_0=0$ or $p_0=\pi$.
But this is impossible since $p_0\in]0,\pi[$. %

(iii) Again using the fundamental equation (\ref{D(x,p)}) we get, 
$$
D(\lambda_j(p),p)=0~~\mbox{  for all $p\in(0,\pi)$ and $j=1,\cdots,N$}.
$$
By differentiating this identity with respect to $p$ and using
$
\frac{\partial D}{\partial p}(x,p)=-2a\sin p,
$
we get
 $$
 \lambda_j'(p)=\frac{2a\sin p}{\Delta'(\lambda_j(p))} 
 ~~\mbox{  for all $p\in(0,\pi)$ and $j=1,\cdots,N$}.
$$ 
So the sign of  $\lambda_j'$  on $]0,\pi[$ is that of $a\Delta'$ on  $]\mu_{j-1},\mu_j[$. 
As $a=(-1)^{N-1}a_1a_2\cdots a_N$ and according to the Lemma \ref{lemma3} where the sign of $\Delta'$ has been found, 
we see that $(-1)^{N+j-1}\lambda'_j(p)>0$ for any $p\in]0,\pi[$. 
 %
 %%%%
\begin{example} Let  $N=2$ so that $A=a_{2}E_{2,1} ~,~B=\mbox{tridiag}(\{a_1\},\{b_1,b_2\})$
where $a_1,a_2>0 $ and $b_1,b_2\in\mathbb{R}$. In this case, it is easy to see on one hand that
 $h(0)$ has two simple eigenvalues.  On the other hand,  $h(\pi)$ has a double eigenvalues if, and only if, 
$b_1=b_2$ and $a_1=a_2$,  in such case the double eigenvalue is $\mu=(b_1+b_2)/2=b_1/2$. 

In the second part of this example we consider the inverse problem of this model. More specifically, 
let $\alpha_1<\beta_1\leq\alpha_2<\beta_2$ and find $a_1,a_2>0 $ and $b_1,b_2\in\mathbb{R}$ such that the corresponding $H_0$ satisfies
$
\sigma(H_0)=[\alpha_1,\beta_1]\cup[\alpha_2,\beta_2].
$
First, if there is a solution then, according to Theorem \ref{thm7.1}, 
$\alpha_1,\beta_2$ are the eigenvalues of $h(0)$  and $\beta_1,\alpha_2$ are the eigenvalues of  $h(\pi)$. Thus
$$
\alpha_1+\beta_2=b_1+b_2=\alpha_2+\beta_1~,b_1b_2=(a_1-a_2)^2+\alpha_2\beta_1~ \mbox{ and }~\alpha_1\beta_2=\alpha_2\beta_1 -4a_1a_2
$$
In particular,  the numbers $\alpha_1,\beta_1,\alpha_2,\beta_2$ are not completely arbitrary. In fact,
once the identity $\alpha_1+\beta_2=\alpha_2+\beta_1$ holds we get 
$16a_1a_2=(\beta_2-\alpha_1)^2-(\alpha_2-\beta_1)^2$; and  $b_1$ and $b_2$ exist if, and only if, the discriminant 
$
(\alpha_2-\beta_1)^2-4(a_1-a_2)^2\geq0.
$
Here again we see that the gap is degenerate, i.e. $\alpha_2-\beta_1=0$  if, and only if, $a_1=a_2=(\beta_2-\alpha_1)/{4}$ and $b_1=b_2=\alpha_2+\beta_1$.
All that is known for general periodic Jacobi matrices, see \cite{CH,H,vM}. The goal here is to show the usefulness of the identity (\ref{D(x,p)}) and to contrast with the model considered in the example \ref{inverse1} where the inverse problem has a unique solution.
\end{example}
%%%%
%%
\protect\setcounter{equation}{0}
 %%%%%%%%%%%%%%%%%%%%%%%%%%%%%%%%%%%%%%%%%%%%%%%%%%%%%%%%
 %%%%%%%%%%%%%%%%%%%   combinaison des deux cas %%%%%%%%%%%%%%%%%%%%%%%
 %%%%%%%%%%%%%%%%%%%%%%%%%%%%%%%%%%%%%%%%%%%%%%%%%%%%%%%%%
 %%%%%%%%%%%%%%%%%%%%%%%%%%%%%%%%%%%%%%%%%%%%%%%%%%%%%%%%%
 \section{A case where $A$ is a singular indefinite matrix}\label{singular3}
Remark first that the assertions of the previous Theorem \ref{thm7.0} remain valid if we replace $A$ by 
$
a_NE_{1N}
$. 
Here  we study the combination of the two cases: we modify the last model by putting  non zero coefficients  in the two corners of $A$.   
While  the band/gap structure of the spectrum of $H_0$ may break down our Mourre estimate still holds. 
More precisely, let $N\geq2$ and $a_1,\cdots,a_{N+1}>0$, $b_1,\cdots, b_N\in\mathbb{R}$. 
Let $H_0=H_0(A,B)$  defined by (\ref{Jacobi1}) where
\begin{eqnarray}\label{A,B,combinaison}
A=a_NE_{N1}+a_{N+1}E_{1N}\quad \mbox{and}\quad B=\mbox{trid}(\{a_{i}\}_{1\leq i\leq N-1},\{b_i\}_{1\leq i\leq N})
\end{eqnarray}
Notice that if $N\geq3$ then  $A$ is a singular indefinite matrix, since its spectrum is nothing but $\{0,\pm\sqrt{a_Na_{N+1}}\}$.
We have,
\begin{theorem} \label{thm7.1'}
\begin{enumerate}
\item There exist $\alpha_1<\beta_1,\alpha_2<\beta_2,\cdots,\alpha_N<\beta_N$  such that  the spectrum of $H_0$ is purely absolutely continuous and
$
\sigma(H_0)=\cup_{j=1}^N[\alpha_j,\beta_j],
$
\item If  (\ref{HP1}) holds then  we have a Mourre estimate for $H_0$ 
on $\mathbb{R}\setminus\kappa(H_0)$, i.e.
$
\mu^{\mathbb{A}}(H_0)=\mathbb{R}\setminus\kappa(H_0)
$
($\kappa(H_0)$ is finite and contains $\{\alpha_1,\beta_1,\cdots,\alpha_N,\beta_N\}$).
 In particular,  conclusions of Theorem \ref{rbvcontinuity} hold for $H_0$\end{enumerate}
\end{theorem}
%
%%%%%%%%%%%%%%%%%%%%%%%%%%%%%%%%%%%%%%%%%%%%%%
It should be clear that one may deduce results for  $H=H_0+V$ for a class of compact perturbations $V$, but we will not state them separately.

Again we will  look at the eigenvalues of $h(p)$
 given in this case  by
$$
h(p)=\left(
\begin{array}{ccccc}
  b_1&a_1&0&\cdots&z_p   \\
  a_1&b_2&a_2&0&\cdots\\
  0&\ddots&\ddots&\ddots&\vdots\\
  0&\cdots&a_{N-2}&b_{N-1}&a_{N-1}\\
  \overline{z_p}&\cdots&0&a_{N-1}&  b_N  
\end{array}
\right)\quad{with}\quad z_p=a_Ne^{ip} + a_{N+1}e^{-ip}.
$$
Since $h(-p)$ is the transpose of ${h(p)}$, it  is enough to study the  operators $h(p)$ for $p\in[0,\pi]$. Theorem \ref{thm7.1'} is a consequence of the following:
\begin{theorem}\label{thm7.0'}
The eigenvalues $\lambda_1(p), \cdots,\lambda_N(p),$ of $h(p)$ are  non constant $[0,\pi]$. 
They are all simple for all $p\in[0,\pi]$ except for finite set of $p$'s for which some of them may be double. 
Such set is reduced to $\{0,\pi\}$ if $a_N\not=a_{N+1}$.
 \end{theorem}
 %
%
%%%%%%%%%%%%%%  Les trois Žtapes %%%%%%%%%%%%%%%%%%%%%%%%
%%%%%%%%%%%%%%%%%%%%%%%%%%%%%%%%%%%%%%%%%%%%%%%
The proof  of Theorem \ref{thm7.0'} will be done in few steps that we state in lemmas.
\begin{lemma}\label{lemma1'} 
(i) The eigenvalues of $h(\pi/2)$ are  all simple.\\
 (ii) If $a_N\not=a_{N+1}$ then
for all $p\in]0,\pi[$, the eigenvalues of $h(p)$ are all simple,
 while for $p=0,\pi$ some of them may be double.
 \end{lemma}
\textbf{Proof  }  Simple adaptation of  the proof of Lemma \ref{lemma1}.
\begin{lemma}\label{lemma2'} Let $a'=(-1)^{N-1}a_1a_2\cdots a_{N-1}(a_N +a_{N+1})$. The characteristic polynomial of $h(p)$ is given by,
\begin{equation}\label{D(x,p)'}
D(x,p)=\Delta(x)-2a_N a_{N+1}\det(\dot{B}-x)\cos^2p +2a'\cos p,
\end{equation}
where, with $\dot{B}$ as in Lemma \ref{lemma2},
\begin{equation}\label{Delta(x)'}
\Delta(x)=\det(B-x)-(a_N - a_{N+1})^2\det(\dot{B}-x)
\end{equation}
\end{lemma}
\textbf{Proof  } Direct calculations.
\begin{lemma}\label{lemma3'}
 The derivative $\Delta'$ of $\Delta$ has $N-1$ roots
 $\mu_1<\mu_2<\cdots<\mu_{N-1}$ and 
the $\lambda_j$'s may be chosen so that, for all $j=1,\cdots,N$,
 $
 (-1)^j\Delta'(\lambda_j(\pi/2))>0.
 $
\end{lemma}
\textbf{Proof }  Again remark that $\Delta(x)=D(x,\pi/2)$ and denote by $\lambda_1(\pi/2)<\lambda_2(\pi/2)<\cdots<\lambda_N(\pi/2)$ the distinct eigenvalues of $h(\pi/2)$. 
 They are roots of $\Delta(x)=D(x,\pi/2)$.  So for any $j=1,\cdots,N-1$ there is
 $\mu_j\in]\lambda_j(\pi/2),\lambda_{j+1}(\pi/2)[$ such that
$
\Delta'(\mu_j)=0,
$
Moreover, since  $\Delta(x)=(-1)^Nx^n+\mbox{lower order}$,  it is clear that $(-1)^j\Delta'(\lambda_j(\pi/2))>0$ for any $j=1,\cdots N$.
 The proof is finished.
\begin{lemma}\label{lemma4'}
There is a neighborhood of $\pi/2$ on which
the eigenvalues $\lambda_j(p)$ of $h(p)$ are monotonic.
\end{lemma}
\textbf{Proof } Let us fix $j=1,\cdots,N$ and put $\mu_0=-\infty, \mu_N=+\infty$.
According to the Lemma \ref{lemma3'}, $\lambda_j(\pi/2)\in]\mu_{j-1},\mu_{j}[$
 and $(-1)^j\Delta'(\lambda_j(p))>0$. We also have 
 $
\frac{\partial D}{\partial p}(\lambda_j(\pi/2),\pi/2)=-2a'(a_N^2+a_{N+1})^2\not=0,
$
Now after differentiation of the identity $D(\lambda_j(p),p)=0$  at $p=\pi/2$ we find,
 $$
 \lambda_j'(\pi/2)=\frac{2a'(a_N^2+a_{N+1}^2)}{\Delta'(\lambda_j(\pi/2))}\not=0.
$$ 
The proof is complete.
%
%%%%
%%%%
%%
%%%%%%%%%%%%%%%%%
%%%%
%%
%
\begin{example} \label{ex3''}
 Let $A=a_{3}E_{3,1}+ a_{4}E_{1,3}~,~B=\mbox{tridiag}(\{a_1,a_2\},\{b_i=0\})$, with  $a_1,a_2,a_3+a_4>0$  and $a_3,a_4\geq0$. 

\textbf{(i)} If $a_3\not=a_4$ then $h(p)$ has only simple eigenvalues for all $p\in(0,\pi)$. 
Moreover,  $h(0)$ has a double eigenvalue
 if, and only if, $a_1=a_2=a_4+a_3$, which is equivalent to $h(\pi)$ has a double eigenvalue. In such case, 
 $-(a_4+a_3)$ (resp. $a_4+a_3$) is a double eigenvalues of $h(0)$ (resp. $h(\pi)$).

\textbf{(ii)} Assume that $a_3=a_4$. In this case, $h(p)=(2a\cos{p})A+B$. Moreover, 
$$
\lambda_j'(p)=-2a\sin{p}\langle AW_j(p),W_j(p)\rangle=-4a_3x_jz_j\sin{p}
$$
where $(x_j,y_j,z_j)$ are the coordinate of $W_j(p)$. This derivative vanishes if, and only if, $p=0,\pi$ or $x_jz_j=0$.
By writing the eigenvalue problem for $h(p)$ we find out that, $x_jz_j=0$ if, and only if, \\
(i) $x_j=0$  leading to $\lambda_j=\pm a_2$, $\pm\cos p_\pm=\mp a_1/2a_3$, provided $a_1\leq2a_3$ or; 
\\
(ii) $z_j=0$ leading to $\lambda_j=\pm a_1$, $\pm\cos p'_\pm=\mp a_2/2a_3$, provided $a_2\leq2a_3$.\\
Moreover, $h(p)$ has a double eigenvalue if, and only if, $a_1=a_2<2a_3$ and  $p=p_\pm=p'_\pm$.

Let us arrange the eigenvalues of $h(p)$ so that 
$\lambda_1(\pi/2)=-\sqrt{a_1^2+a_2^2},\lambda_2(\pi/2)=0,\lambda_3(\pi/2)=\sqrt{a_1^2+a_2^2}$. Using the preceding formula of $\lambda'_j$ one may prove that 
$ \lambda_1'(\pi/2)<0$, $\lambda_2'(\pi/2)>0$ and $ \lambda_3'(\pi/2)<0$.  In particular,
\begin{enumerate}
\item if $a_1>2a_3$ and $a_2>2a_3$ then, for all $j$, $\lambda_j'(p)$ never vanishes on $(0,\pi)$ so they have the sign of  
$\lambda_j'(\pi/2)$. Hence, the spectrum of $H_0$ is absolutely continuous and
$
\sigma(H_0)=[\lambda_1(\pi),\lambda_1(0)]\cup[\lambda_2(0),\lambda_2(\pi)]\cup[\lambda_3(\pi),\lambda_3(0)]
$.
Moreover,  $\lambda_1(0)<\lambda_2(0)$ otherwise one may find $p_0\in[0,\pi/2]$ such that $\lambda_1(p_0)=\lambda_2(p_0)$ 
which contradicts the simplicity of the eigenvalues. Similarly, $\lambda_2(\pi)<\lambda_3(\pi)$ so that we have two gaps.
\item If $a_1\leq 2a_3<a_2$ then there is two points $p_\pm$ given by $\cos p_\pm=\mp a_1/2a_3$ for which  $h(p_\pm)$ has $\pm a_2$ 
as a simple eigenvalue whose corresponding eigenvectors satisfy $x_jz_j=0$. 
  Hence $\kappa(H_0)=\{\pm a_2,\lambda_j(0),\lambda_j(\pi),j=1,2,3\}$.
  For example, if $a_1=a_3=1$ and $a_2=2$ one may deduce that
 the spectrum of $H_0$ is purely absolutely continuous with two gaps and, for some $\alpha>1$,
$
\sigma(H_0)=[-\alpha,-1]\cup[-1/2,1/2]\cup[1,\alpha].
$
\item Finally if $a_1=a_2< 2a_3$ then,there is two points $p_\pm$ given by $\cos p_\pm=\mp a/2a_3$ for which  $h(p_\pm)$ has $\pm a$ as a double eigenvalue whose corresponding eigenvectors satisfy $x_jz_j=0$.  Hence one may deduce that  we only have  overlapping spectral bands. 

\end{enumerate}

\end{example}
\section{ An example where $A$ is a lower triangular matrix: Symmetric fourth-order difference operators}\label{triangulaire}%  
 Let $a_n,b_n>0$ and $v_n\in\mathbb{R}$. Consider the difference operator $\mathcal{J}=\mathcal{J}(a_n,b_n,v_n)$
 acting in $l^2(\mathbb{Z})$ by 
$$
(\mathcal{J}x)_n=a_nx_{n+2}+b_nx_{n+1}+v_nx_n+b_{n-1}x_{n-1}+a_{n-2}x_{n-2}.
$$
Assume that, for some  $a,b,c,d>0$ and  $\sigma\geq0$, one has
$$
\lim_{|n|\rightarrow\infty}|n|^{\sigma}(|a_{2n}-a|+|a_{2n+1}-b|+|b_{2n+1}-c|)+|b_{2n}-d|+|v_{n}|)=0,
$$
Then $\mathcal{J}$ is a bounded symmetric operator. Moreover,
We have,
\begin{theorem}\label{fourth-order}
(i) There exist $\alpha_1<\beta_1$ and $\alpha_2<\beta_2$ such that $\sigma_{ess}(\mathcal{J})=[\alpha_1,\beta_1]\cup[\alpha_2,\beta_2]$.\\
(ii) If $\sigma\geq1$, 
then there is  a self-adjoint operator $\mathbb{A}'$ such that $\tilde{\mu}^{\mathbb{A}'}(\mathcal{J})=\mathbb{R}\setminus\kappa(\mathcal{J})$ 
for some finite set $\kappa(\mathcal{J})$. In particular,  all possible eigenvalues of $\mathcal{J}$ are finitely degenerate and cannot accumulate outside 
$\kappa(\mathcal{J})$.\\
(iii)
Finally,  if $\sigma>1$ then $\mathcal{J}$ has no singular continuous spectrum and conclusions of Theorems \ref{theorem1} and
\ref{theorem11} hold for $\mathcal{J}$.
\end{theorem}
\textbf{Proof }  Consider the unitary operator $U:l^2(\mathbb{Z})\rightarrow l^2(\mathbb{Z},\mathbb{C}^2)$ defined by 
$
(Ux)_n=(x_{2n}, x_{2n+1})
$.
One may verify that $U\mathcal{J}U^{-1}=H(A_n,B_n)$ with
$
A_n=\left(
\begin{array}{cc}
  a_{2n}&0     \\
 b_{2n+1}&  a_{2n+1}
\end{array}
\right) \mbox{,}B=\left(
\begin{array}{cc}
  v_{2n}&b_{2n}   \\
 b_{2n}&  v_{2n+1} 
\end{array}
\right)
$.
Our hypothesis means  that $H=H(A_n,B_n)$ is a compact perturbation of 
 $H_0=H_0(A,B)$   defined by (\ref{Jacobi0}) with 
$
A=\left(
\begin{array}{cc}
  a&0     \\
  c&  b  
\end{array}
\right) \mbox{,}B=\left(
\begin{array}{cc}
  0&d   \\
  d&  0 
\end{array}
\right)$. 
 In this case, 
$
h(p)=\left(
\begin{array}{cc}
  2a\cos{p}&z_p    \\
  \overline{z_p}&  2b\cos{p}  
\end{array}
\right), \mbox{ with } z_p=d+ce^{ip}.
$
One may easily prove that,  $h(p)$ has a double eigenvalue if, and only if, $p=\pi, a=b$ and $c=d$.
Moreover, for any $p\in(0,\pi)$, $h(p)$ has two simple eigenvalues given by 
$$
\lambda_j(p)=(a+b)\cos{p}+(-1)^j \sqrt{(a-b)^2\cos^2{p}+4dc\cos^2{p/2+(d-c)^2}}, j=1,2.
$$
Hence our general method applies directly here also. The proof is complete.
\begin{example}
Suppose that  $a=b$ and  $c=d$. Then  $\lambda_2$ is decreasing and generates the spectral 
band $\Sigma_2=[-2a,2(a+d)]$. In contrast,  we have:
\begin{enumerate}
\item if $K:=\frac{d}{4a}>1$ then $\lambda_1$ is increasing and  $\Sigma_1=[2(a-d),-2a]$.
 In this case, $\kappa(J)=\partial\Sigma_1\cup\partial\Sigma_1$ where  $\partial\Sigma_i$ consists of the edges of  $\Sigma_i$; and 
$\sigma(J_0)=[2(a-d),-2a]\cup[-2a,2(a+d)]=[2(a-d),2(a+d)]$ (the spectral bands touch in one point $-2a$ which is a double eigenvalue of $h(\pi)$). 
\item if $K\leq1$  then $\lambda_1$ has a critical point at $p_0$ given by $\cos{p_0/2}=K$. 
In this case, $\lambda_1(p_0)=-2a(1+2K^2)$ and there is two possibilities.
If $2a\leq d\leq4a$ then  $\Sigma_1=[-2a(1+2K^2),-2a]$ 
so that the two bands $\Sigma_1$ and $\Sigma_2$  touch in one point. In contrary,  if  $0< d\leq2a$ then  
$\Sigma_1=[-2a(1+2K^2),2(a-d)]$ so that  $\Sigma_1\cap \Sigma_2=[-2a,2(a-d)]$. 
 In both cases $\kappa(J)=\partial\Sigma_1\cup\partial\Sigma_1\cup\{\lambda_1(p_0)\}$ and
$\sigma(J_0)=[-2a(1+2K^2),2(a+d)]$.
\end{enumerate}
Of course in this case $\mathcal{J}$ is a compact perturbation of $\mathcal{J}_0=\mathcal{J}(a_n=a,b_n=d,v_n=0)$ that one may study directly. 
Indeed, by Fourier transform, $\mathcal{J}_0$ is unitarily equivalent to the multiplication operator by the function 
$f(p)=2d\cos{p}+2a\cos{2p}$ acting in $L^2(-\pi,\pi)$. 
\end{example}

\begin{example}
Suppose that  $a=b$ and $d\not=c$.  Then 
here again $\lambda_2$ is decreasing and  $\Sigma_2=[-2a+|d-c|,2a+d+c]$ and
\begin{enumerate}
\item if $a<\frac{dc}{2(d+c)}:=K_1$ then $\lambda_1$ is increasing and  $\Sigma_1=[2a-(d+c),-2a-|d-c|]$.
 In this case, $\kappa(J)=\partial\Sigma_1\cup\partial\Sigma_1$ and 
$\sigma(H_0)=[2a-(d+c),-2a-|d-c|]\cup[-2a+|d-c|,2a+d+c]$ (a spectral gap of length $2|d-c|$). 
\item if $a>\frac{dc}{2|d-c|}:=K_2$ then $\lambda_1$ is decreasing and  $\Sigma_1=[-2a-|d-c|,2a-(d+c)]$.
 In this case, $\kappa(J)=\partial\Sigma_1\cup\partial\Sigma_1$ and 
$\sigma(H_0)=[-2a-|d-c|,2a-(d+c)]\cup[-2a+|d-c|,2a+d+c]$.  One may have a non trivial gap in some cases (e.g. $d=10,c=1$ and $a=1$) 
and overlapping bands in others (e.g. $d=2, c=1$ and $a=3/2$).
\item If $K_1\leq a\leq K_2$ then $\lambda_1$ has a critical point at $p_0$ given by $\cos{p_0}=\frac{dc}{8a^2}-\frac{d^2+c^2}{2dc}$. 
It generates the spectral band $\Sigma_1=[\alpha,\beta]$
where $\alpha=\min(\lambda_1(0),\lambda_1(p_0),\lambda_1(\pi))$ and $\beta=\max(\lambda_1(0),\lambda_1(p_0),\lambda_1(\pi))$.
 In this case,  $\kappa(J)=\partial\Sigma_1\cup\partial\Sigma_1\cup\{-2a-4K^2\}$ and 
$\sigma(H_0)=[\alpha,\beta]\cup[-2a,2(a+d)]$.  For example, let $c=2,d=1$. Then $K_1\leq a\leq K_2$ means that $1/3\leq a\leq 1$. 
$\Sigma_1=[-9/4,-2]$ and $\Sigma_2=[-1,3]$ and we have a spectral gap.
$\Sigma_1=[-3,-1]$ and $\Sigma_2=[-2,4]$ and we have overlapping  spectral bands.
\end{enumerate} 
\end{example}
\begin{example}
It remains to study the case where $a\not=b$. To be simple we only discuss the following two situations. 
\begin{enumerate}
\item If $c=d=0$ then $\lambda_1(p)=2a\cos{p}$ and $\lambda_2(p)=2b\cos{p}$ so that we have two spectral bands one contained in the another.
\item 
In contrast if $d=0$ and $c>0$ then  both $\lambda_j's$ are monotonic and
one may choose $a,b,c$ to get overlapping bands or  bands with a spectral gap between them.
Here $\kappa(J)$ is  the set of the edges of the spectral bands.
\end{enumerate}
\end{example}
\section{Further developments}\label{further}
In this section we explain how our analysis extends to compact perturbation of difference operators of higher order with matrix coefficients.
More specifically, let  
$
H=H_0+V
$
where $V$  is a compact operator in ${\mathcal H}=l^2(\mathbb{Z}, \mathbb{C}^N) $ and 
$H_0$ is defined   on $\mathcal{H}$ by
\begin{equation}\label{Jacobi0'}
(H_0\psi)_n=B\psi_n+\sum_{k=1}^mA_k\psi_{n+k}+A^*_k\psi_{n-k}, ~~\mbox{
for all } n\in\mathbb{Z}.
\end{equation}
where $A_1,\cdots,A_k$ and $B=B^*$ are $N\times N$ matrices.
By using Fourier transform we easily show  that $H_0$ is unitarily equivalent  to the direct integral 
\begin{equation}\label{direct'}
\hat{H_0}=\int_{[-\pi,\pi]}^\oplus h(p)dp\quad\mbox{ acting in }\quad
\hat{\mathcal{H}}=\int_{[-\pi,\pi]}^\oplus\mathbb{C}^Ndp
\end{equation}
where the reduced operators $h(p)$ are given by,
$$
h(p)=B+\sum_{k=1}^me^{-ikp}A_k+e^{ikp}A^*_k.
$$
So by repeating  verbatim the arguments of Sections \ref{libre} and \ref{mestimateH0} we prove 
that Corollary \ref{moureH_0}
and it consequence Theorem \ref{rbvcontinuity}  are valid for $H_0$. Similarly, the perturbative argument developed in Section \ref{perturbŽ} extends to $H=H_0+V$ with similar 
assumptions on the perturbation $V$.
In particular, analogous results to Theorem \ref{MEFORH}  and those stated in the introduction can be obtained.

We finish by given the simplest but non trivial example one may consider here. Assume that $A_k=A$ for all $k=1,\cdots,m$ where $A$ is a positive definite matrix. Then 
$h(p)=B+(D_m(p)-1)A$ where $D_m(p)=1+2\sum_{k=1}^m\cos(kp)$ is the Dirichlet Kernel. So by repeating the argument  of Lemma \ref{deŽrivŽee} we get that
$$
\lambda_j'(p)=D_m'(p)<AW_j(p),W_j(p)>
$$
The Dirichlet kernel is known to strongly oscillates, especially when $m$ is large, so is for $\lambda_j(p)$.
 Hence $\kappa(H_0)$ will be little bit complicate but finite and our general framework applies directly, compare however with section \ref{positivedefinie}.
%%%%%%%%%%%%%%%%%%%%%%%%%%%%%%%%%%%%%%%%%%%%%
%%%%%%%%%%%%%%%%%%%%%%%%%%%%%%%%%%%%%%%%%%%%%
%%%%%%%%%%%%%%%%%%%%%%%%%%%%%%%%%%%%%%%%%%%%%%%%%%%%%%%%%%%%%%%%
%%%%%%%%%%%%%%%%%%%%%%%%%%%%%%%%%%%%%%%%%%%%%%%%%%%%%%%%%%%%%%
%%%%%%%%%%%%%%%%%%%%%                   %%%%%%%%%%%%%%%%%%%%%%
%%%%%%%%%%%%%%%%%%%%%     References    %%%%%%%%%%%%%%%%%%%%%%
%%%%%%%%%%%%%%%%%%%%%                   %%%%%%%%%%%%%%%%%%%%%%
%%%%%%%%%%%%%%%%%%%%%%%%%%%%%%%%%%%%%%%%%%%%%%%%%%%%%%%%%%%%%%

\end{document}